\documentclass[11pt]{article}
\usepackage[mathscr]{eucal}
\usepackage{xypic}
\xyoption{all}
\usepackage{amsfonts}
\usepackage{amsmath}
\usepackage{amsthm}
\usepackage{amssymb}
\usepackage{latexsym}
\usepackage{eufrak}
\newcommand{\fix}{{\mathtt{fix}}}
\newcommand{\TT}{{\mathbb{T}}}
\newcommand{\cht}{{\chi_{_\TT}}}

\newcommand{\dd}{{\mathcal{D}}}
\newcommand{\m}{{\mathfrak{m}}}

\newcommand{\wfrac}{\mbox{$\frac{1}{|W|}$}}
\newcommand{\perv}{{\mathcal{P}\textit{erv}}}

\newtheorem{theorem}[equation]{Theorem}
\newtheorem{proposition}[equation]{Proposition}
\newtheorem{corollary}[equation]{Corollary}
\newtheorem{lemma}[equation]{Lemma}
\theoremstyle{definition}                               

\newtheorem{remark}[equation]{Remark}
\newtheorem{conjecture}[equation]{Conjecture}

\newcommand{\dis}{{\displaystyle}}

\newcommand{\beq}{\begin{equation}\label}

\newcommand{\iso}{{\;\;\stackrel{_\sim}{\longrightarrow}\;\;}}
\newcommand{\cd}{\!\cdot\!}

\newcommand{\vi}{${\sf {(i)}}\;$}
\newcommand{\vii}{${\sf {(ii)}}\;$}
\newcommand{\viii}{${\sf {(iii)}}\;$}

\newcommand{\mr}{|_{\hreg}}

\newcommand{\sset}{\subset}
\newcommand{\ssminus}{\smallsetminus}
\newcommand{\G}{\Gamma}

\newcommand{\g}{{\mathfrak{g}}}
\newcommand{\bg}{{\boldsymbol{{\mathfrak{g}}}}}
\newcommand{\tg}{{\boldsymbol{\widetilde{\boldsymbol{\mathfrak{g}}}}}}
\newcommand{\GG}{{\mathbf{G}}}
\newcommand{\zz}{{\mathcal{Z}}}
\newcommand{\bb}{{\mathcal{B}}}
\newcommand{\Gr}{{\mathcal{G}}r}

\newcommand{\spher}{^{^{\mathtt{spher}}}}

\newcommand{\irrep}{{\sf{Irrep}}}

\newcommand{\id}{{{\mathtt {Id}}}}
\newcommand{\kz}{{{\mathtt {KZ}}}}

\newcommand{\too}{\,\,\longrightarrow\,\,}
\newcommand{\onto}{\,\,\twoheadrightarrow\,\,}

\newcommand{\ad}{{\mathtt{{ad}}^{\,}}}

\newcommand{\Spec}{{\mathtt{Spec}}}

\newcommand{\Hom}{{\mathtt{Hom}}}

\newcommand{\grd}{{\mathtt{gr}}}
\newcommand{\Ext}{{\mathtt{Ext}}}

\newcommand{\Tr}{{{\mathsf {Tr}}}}

\newcommand{\hh}{{\mathsf{H}}}

\newcommand{\ehe}{{\mathbf{e}\mathsf{H}_c\mathbf{e}}}

\newcommand{\e}{{\mathbf{e}}}
\newcommand{\emi}{{\mathbf{e}_{_{^{\boldsymbol{-}}}}}}

\newcommand{\sym}{{\mathsf{Sym}}}

\newcommand{\h}{{\mathfrak{h}}}

\newcommand{\HH}{{\mathcal{H}}}

\newcommand{\A}{{\mathcal{A}}}

\newcommand{\bh}{{\mathbf{h}}}

\newcommand{\hreg}{{\mathfrak{h}^{^{_{\mathsf{reg}}}}}}

\newcommand{\triv}{{\mathsf{triv}}}
\newcommand{\sign}{{\mathsf{sign}}}

\def\C{{\mathbb{C}}}

\def\rr{{\mathscr{R}}}
\def\Z{{\mathbb{Z}}}

\def\oo{{\mathscr O}}

\def\bb{{\mathcal{B}}}

\def\g{{\mathfrak{g}}}

\def\sll2{{\mathfrak{s}\mathfrak{l}}_2}

\def\OH{{\mathcal{O}}_{_{{\sf{Hilb}}^n(\C^2)}}}

\def\ccirc{{{}_{\,^{^\circ}}}}

\makeatletter
\def\downroundfill{$\m@th \setbox\z@\hbox{$\braceld$}%
  \braceld\leaders\vrule height\ht\z@ depth\z@\hfill\bracerd$}
\def\uproundfill{$\m@th \setbox\z@\hbox{$\braceld$}%
 \bracelu\leaders\vrule height\ht\z@ depth\z@\hfill\braceru$}
\def\overround#1{\mathop{\vbox{\m@th\ialign{##\crcr\noalign{\kern3\p@}
      \downroundfill\crcr\noalign{\kern3\p@\nointerlineskip}
      $\hfil\displaystyle{#1}\hfil$\crcr}}}\limits}
\def\underround#1{\mathop{\vtop{\m@th\ialign{##\crcr
      $\hfil\displaystyle{#1}\hfil$\crcr\noalign{\kern3\p@\nointerlineskip}
      \uproundfill\crcr\noalign{\kern3\p@}}}}\limits}
\makeatother
\begin{document}
\setlength{\parindent}{6mm}
\setlength{\parskip}{3pt plus 5pt minus 0pt}
\centerline{\Large {\bf
Finite dimensional representations of  rational}}
\vskip 2pt
\centerline{\Large {\bf Cherednik algebras}}
\vskip 4mm
\centerline{\large {\sc {Yuri Berest, Pavel Etingof and Victor Ginzburg}}}

\vskip 2pt
\begin{abstract} A complete classification and  character
formulas for finite-dimensional irreducible representations of
the rational Cherednik algebra of type $\mathbf{A}$ are given.
Less complete results for other types are obtained. Links to the geometry
of affine flag manifolds and Hilbert schemes are discussed.
\end{abstract}

\section{Main results}

\subsection{Preliminaries}
Fix a finite Coxeter group $W$ in a
complex vector space $\h$. Thus, $\h$ is the complexification
of a  real Euclidean  vector space and $W$ is
generated by a finite set $S\subset W$ of
reflections $s\in S$ with respect to certain
hyperplanes $\{H_s\}_{s\in S}$
in that  Euclidean  space.

For each $s\in S$, we choose 
a nonzero
linear function $\alpha_{s}\in\h^*$ which vanishes on $H_s$
(the positive root corresponding to $s$),
and let $\alpha_s^{\vee}=2(\alpha_s, -)/(\alpha_s,\alpha_s) \in \h $
be the corresponding coroot.
The group $W$ acts  naturally  on the set $S$ by conjugation. 

Put $\ell:=\dim_{_\C}\h$.
Let $\wedge^i\h$, $i=0,1,...,\ell,$
denote the $i$-th wedge power of the tautological reflection
 representation.
Thus, $\wedge^0\h=\triv$ is the trivial 1-dimensional
 representation of $W$, and  $\wedge^{\ell}\h=\sign$
 is the sign representation. It is known that if
 $\h$ is an irreducible $W$-module (as we will assume), then
the representations $\wedge^i\h$ 
are irreducible for all $i$.  

According to [EG], for each
 $W$-invariant function $c: S \to \C\,,\, c\mapsto c_s,$
 one defines a {\it rational Cherednik algebra} $\hh_c=\hh_c(W)$ to be
an associative algebra
generated
by the vector spaces $\h$, $\h^*,$ and the set
$W,$ with defining relations  
(cf.  formula \cite[(1.15)]{EG} for $t=1$) given
by
\begin{equation}
\label{defrel}
\begin{array}{lll}\displaystyle
&{}_{_{\vphantom{x}}}w\cd x\cd w^{-1}= w(x)\;\;,\;\;
w\cd y\cd w^{-1}= w(y)\,,&
\forall y\in \h\,,\,x\in \h^*\,,\,w\in W\break\medskip\\
&{}^{^{\vphantom{x}}}{}_{_{\vphantom{x}}}x_1\cd x_2 = 
 x_2\cd x_1\;\enspace,\enspace\;
y_1\cd y_2=y_2\cd y_1\,, &
\forall y_1,y_2\in \h,\;x_1\,,\,x_2 \in \h^*\,\break\medskip\\
&{}^{^{\vphantom{x}}}y\cd x-x\cd y = \langle y,x\rangle
-\!\!\!\sum\limits_{^{_{s\in S}}}
c_s\cd\langle y,\alpha_s\rangle
\langle\alpha_s^\vee,x\rangle \cd s\,,& \forall y\in
\h\,,\,x\in \h^*\,.
\end{array}
\end{equation}
 Thus, the elements  $x\in\h^*$ generate a 
subalgebra $\C[\h]\subset \hh_c$ of polynomial functions on $\h$,
 the elements  $y\in\h$ generate a 
subalgebra $\C[\h^*]\subset \hh_c$,
and the elements $w\in W$
span  a copy of the group
algebra $\C W$ sitting naturally inside $\hh_c$.

Given a finite dimensional  $W$-module $\tau$, extend it
to a module $\tilde\tau$ over $\C[\h^*]\#W$, the crossed product algebra,
by letting $\C[\h^*]$ act via the evaluation map $P\mapsto P(0)$.
Following [DO] and [BEG],
define a standard $\hh_c$-module $M(\tau)$ to be the induced
module $M(\tau):= \hh_c\otimes_{_{\C[\h^*]\#W}}\,{\tilde\tau}.$ This module is
known, by [DO], to have a unique simple quotient, to be denoted $L(\tau)$.
Furthermore, in [OR] 
an abelian category  $\oo(\hh_c)$
(which was denoted by $\oo_0(\hh_c)$ in [BEG]) of `bounded below'
$\hh_c$-modules  has been introduced,
 similar to one defined by Bernstein-Gelfand-Gelfand
in the case of semisimple Lie algebras. By definition, the objects of
  $\oo(\hh_c)$ are finitely-generated $\hh_c$-modules $M$ such that
the $\h$-action on $M$ is locally-nilpotent.
The modules $\{L(\tau)\}_{\tau\in \irrep(W)}$
were shown to be precisely the simple objects of the category  $\oo(\hh_c)$.

In this paper we will be interested in classification and characters 
of {\it finite-dimensional}
$\hh_c$-modules. They clearly belong to the category
$\oo(\hh_c)$, but exist only for a certain special
set of parameters `$c$'. 

\begin{remark} Let $\varepsilon: W\to \Bbb C^\times$ be a group
character. 
Then there exists an isomorphism $\hh_c\iso \hh_{\varepsilon\cdot  c}$
acting identically on $\h$ and $\h^*$, and sending each $w\in W$ to $\varepsilon(w)\cdot w$,
so the representation categories of $\hh_c$ and $\hh_{\varepsilon\cdot  c}$ are equivalent. 
Under this equivalence, representations with lowest weight $\tau$ go to representations 
with lowest weight $\varepsilon\otimes\tau$. That explains symmetry
between our results for $\hh_c$ and for $\hh_{\varepsilon\cdot  c}$,
at various places below. $\lozenge$ 
\end{remark} 
\subsection{Type $\mathbf{A}$.}
We start with the case 
$W=S_n$, the Symmetric group, and $\h=\Bbb C^{n-1}$, 
so that the  parameter `$c$' runs over the line $\C$.

The following result provides a complete classification 
of irreducible finite dimensional representations of the algebra
$\hh_c=\hh_c(S_n)$.

\begin{theorem}[Classification  of finite dimensional
representations in type ${\mathbf{A_{n-1}}}$]\label{classif}\hfill\break

\vskip -10mm
\vi The only values of $c\in \C$ for which nonzero 
finite dimensional representations of $\hh_c$ exist are the rational
numbers of the form:
$\pm r/n$, where  $r\in \Bbb N$ with $(r,n)=1$.

\vii Let $c=\pm r/n$ with $r\in \Bbb N$, $(r,n)=1$. Then the 
only irreducible finite dimensional representation of $\hh_c$ is
$L(\triv)$, in the `$+$' case, resp. $L(\sign)$, in  the `$-$' case.
\end{theorem}
%\viii The character of $L(\varepsilon)$ is given by the formula 
%$$
%\chi_{L(\varepsilon)}(t,w)\;=\;\varepsilon(w)
%\cdot t^{-(r-1)(n-1)/2}\cdot\frac{\det(1-t^r\cd
%w)}{\det(1-t\cd w)} \quad,\quad\forall w\in W\,.
%$$
%\end{theorem}

\begin{remark}
Some special cases of this Theorem can be found 
in [CO], [Dz], and [Go]. Also, it was pointed out by Cherednik 
\cite[p.65]{Ch2} that 
Theorem \ref{classif} can be deduced from the 
results of [Ch2] (see Section \ref{Chere} below for more details). $\lozenge$
\end{remark}
\subsection{Geometric realization}
Let ${\mathbb P}^k$ denote the complex projective space
of dimension $k$, and let $\C_{_{{\mathbb P}^k}}$ denote the constant
sheaf on ${\mathbb P}^k$ (extended by zero, whenever
viewed as a sheaf on ${\mathbb P}^l\supset{\mathbb P}^k$). 
Let $\perv({\mathbb P}^{n-1})$ be the abelian  
 category of perverse sheaves on ${\mathbb P}^{n-1}$ which are
constructible with respect to the standard stratification
${\mathbb P}^{n-1}= {\sf{pt}}\,\cup\, \C^1\,\cup\, \C^2\,\cup\ldots\,
\cup\,\C^{n-1}.\,$ The perverse sheaves
$\C_{_{{\mathbb P}^k}}[k]\,,\, k=0,1,\ldots,n-1,$
are exactly the simple objects of the category  $\perv({\mathbb P}^{n-1})$.

\begin{theorem}\label{perv} Let $W=S_n$, and $c
=r/n$ with $(r,n)=1\,,\,r\in\Bbb N$. Then 

\vi The block of the category $\oo(\hh_c)$ containing the
(only) finite dimensional representation is equivalent
to the abelian  
 category $\perv({\mathbb P}^{n-1})$.
The equivalence  sends the perverse sheaf
$\C_{_{{\mathbb P}^k}}[k]$ to
the $\hh_c$-module $L(\wedge^k\h)\,,\,\forall k=0,1,\ldots,n-1.$

\vii The class  of $L(\wedge^k\h)$ in the Grothendieck group of the
category $\oo(\hh_c)$
is given by the formula 
$[L(\wedge^k\h)]=\sum_{j=k}^n\,(-1)^{j-k}\cdot [M(\wedge^j\h)]$,
for any $k=0,1,\ldots,n-1.$

\viii Any  block in  $\oo(\hh_c)$ other than the one considered
in \vi is  semisimple.
\end{theorem}

\begin{remark} If $c=-r/n$, the category $\oo(\hh_c)$ has a similar structure,
since we have an isomorphism $\hh_c\simeq\hh_{-c}$. $\lozenge$
\end{remark}

\subsection{Other types} 

We are going to generalize Theorem \ref{classif} to other Coxeter groups $W$. 
Given a group homomorphism $\varepsilon: W\to\C^\times$, introduce the following

\noindent
{\bf Notation:\;} $\;r:=\frac{2}{\ell}\cdot\sum_{s\in S}\, c_s\cdot\varepsilon(s)$. 

\noindent
For example, if $c$ and $\varepsilon$ are constant functions
on $S$, then $c=\varepsilon\cdot  r/h$, where we
put $h:=2\cdot|S|/\ell$. It is well known, see e.g. \cite[3.18]{hum}, that the
number $h$ thus defined is equal to the Coxeter number of $W$. 

\begin{theorem} \label{BD1}
Assume that 
$c,\varepsilon$ are constant on $S$, and we are in one of the following situations. 

\vi\; $\;r=k\cdot h+1$, 
where $k$ is a nonnegative integer; 

\vii\; $W=W({\mathbf{B}}_n)$, $r=2k+1$, where $k$ is a nonegative integer, 
such that $(r,n)=1$.

\viii\; $W=W({\mathbf{D}}_n)$, $r=2k+1$, where $k$ is a nonegative integer,
such that $(r,n-1)=1$. 

\noindent
Then the only finite dimensional irreducible $\hh_c$-module is $L(\varepsilon)$.
\end{theorem}

The proof of (a generalization of) this theorem  
is given in Section \ref{BDres}.
For analogues of Theorem \ref{BD1}
for the double-affine Hecke algebra see [Ch2, \S8].

\begin{remark} For types other than type ${\mathbf{A}}$ 
there may exist more than one irreducible finite dimensional representation
of $\hh_c$ for a given $c$, and they may have lowest weights
which are not 1-dimensional. As we will see below, 
it happens for example, in type ${\mathbf{D}}_4$. This
seems to be an `$L$-packet' phenomenon, cf. \S7.1. $\lozenge$
\end{remark}
\smallskip

Let $\{x_i\}$ and $\{y_i\}$ be dual bases of $\h$ and $\h^*$,
respectively. Following [BEG], put
$\bold h:=\frac{1}{2}(\sum_{i=1}^\ell\, x_iy_i+y_ix_i)\in \hh_c$.
This element commutes with $W$ in $\hh_c$.
It is known that the $\bh$-action on any  $V\in \oo(\hh_c)$
is locally finite, moreover, for each $a\in \C$,
the generalized $\bh$-eigenspace $V_a$ corresponding to the eigenvalue $a$
is finite dimensional (and, of course, $W$-stable). These (generalized) eigenspaces
$V_a$ will be referred to as {\it weight} spaces of $V$,
so  $V=\oplus_a\, V_a$, where the sum ranges (by [BEG])
over a finite union of (one-sided)
arithmetic progressions with step one each. 

For a finite dimensional $W$-module $E$ and  $w\in W$, write
$\Tr(w, E)$ for the trace of $w$-action in $E$. Now,
given $V\in \oo(\hh_c)$ with weight decomposition  $V=\oplus_a\, V_a$,
and  $w\in W$, define the character of $V$
as a  formal infinite series in (complex) powers of $t$ given by:
$
\Tr_{_V}(w\cdot  t^{\bold h})\,:=\sum_a\, \Tr(w, V_a)\cdot t^a.$
For example, the character of the standard module
$M(\tau)$ is 
equal to
\begin{equation}\label{kapa}
\Tr_{_{M(\tau)}}(w\cdot  t^{\bold h})
=\frac{t^{\kappa(c,\tau)}\cdot\Tr(w,\tau)}{\det(1-t\cdot w)}\,,\quad\text{where}\quad
\kappa(c,\tau)=\frac{\ell}{2}-\sum_{s\in S}\,c_s\cdot s|_{\tau}.
\end{equation}
Here and elsewhere,
`$\det$' denotes the determinant of the reflection representation,
and `$\sum_{s\in S}\,c_s\cdot s|_{\tau}$' denotes the scalar by
which the indicated central element of $\C{W}$ acts in the simple
$W$-module $\tau$. Formula \eqref{kapa} follows from the
natural identification
$M(\tau)\simeq \C[\h]\otimes\tau,$ and the fact that the lowest
weight of $M(\tau)$ equals $\kappa(c,\tau),$ see e.g. [BEG].

\begin{theorem}[Character formula]\label{char}
Under the assumptions of either Theorem \ref{classif} or
Theorem \ref{BD1}\, 
the representation $L(\varepsilon)$ has the character:
$$
\Tr_{L(\varepsilon)}(w\cdot t^\bh)\;=\;\varepsilon(w)
\cdot t^{-(r-1)\cdot\ell/2}\cdot\frac{\det(1-t^r\cd
w)}{\det(1-t\cd w)} \quad,\quad\forall w\in W\,.
$$
In particular, $\dim L(\varepsilon)=r^\ell$.\hfill\qed\break
\end{theorem}

For the rest of the subsection assume that $W$ is a Weyl group of a simple 
Lie algebra.  Write $Q$
for the root lattice in $\h^*$. 
The group  $W$ acts naturally on  $Q$.
This induces, for any integer $r>0$, an $W$-action on
the finite set $Q/r\cd Q$,
making the $\C$-vector space $\C[Q/r\cd Q]$, of $\C$-valued
functions on $Q/r\cd Q$,
a $W$-module. 
The following corollary of the character formula above
is a generalization of \cite[Conjecture 5.10]{BEG}.
 
\begin{proposition}\label{traceCorollary} 
In the situation of either Theorem \ref{classif} or Theorem  \ref{BD1}, 
and $c>0$, one has:
 $L(\varepsilon)|_{W} \,\simeq\,\C[Q/r\cd Q]$.
\end{proposition}

\begin{proof} 
E. Sommers showed in [So] (using a case by case
  argument) that 
if $(r,h)=1$ then,
for any $w\in W$, one has $\Tr(w,\Bbb C[Q/r\cd Q])=r^{\fix(w)}$, where $\fix(w)$
is the dimension of the space of invariants of the
$w$-action in $\h$. 
On the other hand, we calculate
\begin{equation}\label{t=1}
\Tr\bigl(w, L(\varepsilon)\bigr)\,=\,\Tr\bigl(w\cdot t^\bh\,,\,
L(\varepsilon)\bigr)|_{t=1}\; \stackrel{{\tt{Thm.}} \ref{char}}{=\!=\!=}\;
%= \quad\text{(by Theorem \ref{char})}\\ &=\,
\lim_{t\to 1}\,\frac{\det(1-t^r\cdot w)}{\det(1-t\cdot
w)}\,=\,r^{\fix(w)}\,.
%=\,\Tr(w,\Bbb C[Q/r\cd Q]). \end{align*}
\end{equation}
Hence, the $W$-modules $L(\varepsilon)$ and $\Bbb C[Q/r\cd Q]$ have equal
characters. The result follows.
\end{proof}

Given a $W$-module $E$ and $\tau\in\irrep(W)$, let  ${}^{\tau\!\!}{E}:=
(\tau\otimes E)^W$
denote  the $\tau$-isotypic
component of
 $E$.
\begin{corollary}\label{isotypic}
 In the situation of either Theorem \ref{classif} or Theorem  \ref{BD1}, 
and $c>0$,  we have:
$$\dim\bigl({}^{\tau\!\!}{L(\varepsilon)}\bigr) \,=\,
\wfrac\cdot\sum\nolimits_{w\in W}\,\Tr(w,\tau)\cdot 
r^{\fix(w)}\,.
$$
\end{corollary}

\begin{proof} Use \eqref{t=1} and the  general
formula:
$\dim({}^{\tau\!\!}{E})=
\frac{1}{|W|}\cdot\sum_{w\in W}\,\Tr(w,\tau)\cdot\Tr(w,E).$
\end{proof}
\subsection{Representations of the spherical subalgebra}
An  analogue of
Theorem \ref{classif}  also holds
for finite dimensional representations of the {\it spherical subalgebra}
$\ehe\subset \hh_c$, where $\e=\frac{1}{|W|}\sum_{w\in W}\,w,$ see
[EG],[BEG]. 
To formulate it, 
note that, for any $\hh_c$-module
$V$, the space $\e\cd V$ has a natural $\ehe$-module structure,
and the element $\bh$ preserves the subspace $ \e\cd V\subset V$.
Also, write $d_i\,,\,i=1,\ldots,\ell,$ for the
degrees of the basic $W$-invariant polynomials on
$\h$. 

\begin{theorem}[Character formula for $\ehe$]\label{char_ehe} 
\vi In type ${\mathbf{A}}_{n-1}$,  
the only {\it positive} values of $c$ for which nonzero 
finite dimensional representations of $\ehe$ exist are 
$r/n$, where  $r=1,2,\ldots,$ with $(r,n)=1$. 

\vii Under the assumptions of either Theorem \ref{classif} or
Theorem \ref{BD1},\, $\e\cdot L(\triv)$ is the only finite-dimensional
simple $\ehe$-module, if $c>0$.

\viii We have
$$
\Tr_{\e\cdot L(\triv)}(t^{\bh})\,=\,t^{-(r-1)\ell/2}\cdot\prod_{i=1}^\ell
\frac{(1-t^{d_i+r-1})}{(1-t^{d_i})}.
$$
In particular,  we have:
$\dim(\e\cd L(\triv))= \prod_i \frac{d_i+r-1}{d_i}$ (= $\frac{1}{r}\cd 
\left({r+n-1\atop n}\right)$ in  type ${\mathbf{A}}_{n-1}$).
\end{theorem}

\begin{remark} Finite dimensional representations of $\ehe$ for $c<-1$
have the same structure as those for $c>0,$ since the algebra $\e\hh_c\e$ is
isomorphic to $\e\hh_{-c-1}\e$, see e.g. [BEG2]).
For $W=S_n$ and
$-1\le c\le 0,$  there are no finite dimensional representations.
Indeed, if $V$ were a nonzero finite dimensional representation of $\ehe$ then
$V'=\hh_c\e\otimes_{\ehe}V$ would have been
 a nonzero finite dimensional representation
of $\hh_c$ such that $\e\cdot 
 V'\neq 0.$  On the other hand, from Theorem 
\ref{classif} it follows that
such representation does not exist: $\e\cdot L(\sign)=0$ in this case.
$\lozenge$\end{remark} 
\smallskip 

We postpone the proof of parts (i),(ii) of  Theorem \ref{char_ehe}
until later sections.  
However, 
assuming that the character
formula for $L(\triv)$ is already known, part (iii) of the Theorem
is an immediate corollary of
the following general lemma from the theory of complex
reflection groups.  

\begin{lemma}[Solomon \cite{S}]\label{sollem}  Let $W$ be a complex reflection
group acting in a complex vector space $\h$ of dimension $\ell$. Then
one has an equality of rational functions in $t,y$:
$$
\wfrac\sum_{w\in W} \frac{\det(1+y\cdot w)}{\det(1-t\cdot
  w)}=\prod_{i=1}^\ell\frac{1+y\cdot t^{d_i-1}}{1-t^{d_i}}.
$$
\end{lemma}

\noindent
{\sl Proof of Lemma.\;} The LHS is the two-variable Poincar\'e series
of the space $\Omega^\bullet(\h)^W$ of invariant polynomial differential forms on $\h$
(powers of $y$ count the rank of the form and powers of $t$ the
homogeneity degrees of the  polynomial coefficients). However, as was shown in [S],
$\Omega^\bullet(\h)^W$ is a free supercommutative algebra 
generated by basic invariants $p_1,...,p_\ell$ and their
differentials $dp_1,...,dp_\ell$.\footnote{
Recall that $\h/W:=\Spec \C[\h]^W$ is  an affine space with
coordinates $p_1,...,p_\ell$, by Chevalley theorem.
The result of Solomon
may be seen as claiming that the projection $\pi: \h\to\h/W$ induces an isomorphism $\pi^*:
\Omega^\bullet(\h/W)\iso\Omega^\bullet(\h)^W$.
Here is a sketch of proof of this result different from the one given in \cite{S}.

To prove the claim
it suffices to verify that  $\pi^*$ is  an isomorphism 
 on the complement of a codimension 2
subset in $\h/W$.
Further, the claim is clear over the generic point
of  $\h/W$. Hence, one is reduced to proving
the claim at the generic point of a root hyperplane.
The latter case is essentially the one-dimensional situation
that can be easily verified directly. \qed}
As the RHS of the equality in
the lemma is exactly the Poincar\'e series of such free algebra,
the Lemma follows. \qed

To obtain Theorem \ref{char_ehe}(iii), it suffices 
to use the character formula of Theorem \ref{char},
note that $\Tr_{_{\e\cdot L(\triv)}}(t^\bh)=
\frac{1}{|W|}\cdot\sum_{w\in W}\, \Tr_{_{L(\triv)}}(w\cdot t^\bh),$
and to substitute
$y=-t^r$ in Lemma~\ref{sollem}. 

\begin{remark} Compatibility of the dimension formula for
$\e\cdot L(\triv)$ arising from Theorem \ref{char_ehe}(iii)
with the one given, in the special case $\tau=\triv$, by
Corollary \ref{isotypic} is equivalent, modulo the standard identity
$\prod_i\, d_i=|W|$, to the
Shephard-Todd formula:
$\sum\nolimits_{w\in W}\,r^{\fix(w)}= \prod_{i=1}^\ell\,(r+d_i-1).$
$\lozenge$
\end{remark}

\subsection{Absence of self-extensions}
\begin{proposition}\label{ext}
\vi For any simple object $V\in \oo(\hh_c)$ one has:
$\Ext^1_{\oo(\hh_c)}(V,V)=0$. In particular, 

\vii In all
the situations considered in this section, any finite-dimensional 
$\hh_c$-module is a multiple of a (unique) simple $\hh_c$-module.
\end{proposition}

\begin{proof}
Let $\tau\in \irrep(W)$ be the lowest weight of $V$.
Let $N$ be an extension of $V$ by $V$, and  $N_0$
the lowest nonzero generalized eigenspace  of $\bold h$ in $N$.
Any nonzero vector $v\in N_0$ is a 
singular vector, belonging to the 
$W$-module $\tau$, since the $W$-action is semisimple.
Therefore, it generates a proper $\hh_c$-submodule in $N$
(because it is a quotient of $M(\tau)$, it cannot contain 
the whole $N_0$), which has
to be $V$. This shows that $N=V\oplus V$.
\end{proof}

\subsection{The Gorenstein property of $L(\triv)$.}

The action of the  commutative
subalgebra $\C[\h]\sset \hh_c$
on the lowest weight vector in the
simple module $L(\triv)$ clearly
generates the whole space $L(\triv)$. Therefore one has
$L(\triv)=\C[\h]/I$, where $I$ is an ideal in  $\C[\h]$.
Thus, $L(\triv)$ becomes a positively graded commutative algebra.

\begin{proposition}\label{gor}
If $\dim L(\triv) <\infty,$  then $L(\triv)=\C[\h]/I$ is a
Gorenstein algebra.
\footnote{For a geometric interpretation of the Proposition cf. \S7.1.}
\end{proposition}
\begin{proof}
Observe first that the module $L(\triv)$ has a canonical (Jantzen type)
contravariant nondegenerate bilinear form $J(-,-)$, see [DO], [Dz].
To make this form {\it invariant} rather than {\it contravariant},
we use the $SL_2(\C)$-action on any finite dimensional
$\hh_c$-module that has been constructed in \cite[Remark after
Proposition 3.8]{BEG}. In particular, let $F$ be the `Fourier transform'
endomorphism of $L(\triv)$ corresponding to the
action of the matrix $\left(\begin{matrix} \,0 & 1\\ -1 & 0
\end{matrix}\right)\in SL_2(\C)$.
Since the Fourier transform interchanges 
$\h$ and $\h^*$, the bilinear form: $\Phi: (v_1,v_2)\longmapsto J(v_1\,,\,Fv_2)$
on  $L(\triv)$ is compatible with the algebra
structure on $L(\triv)=\C[\h]/I,$
i.e. for any $x\in \h^*$, we have $\Phi(xv_1,v_2)=\Phi(v_1,xv_2)$. 
The $SL_2$-representation theory now shows that
 the top degree homogeneous component
of the graded algebra  $\C[\h]/I$ is 1-dimensional.
Moreover,
the form $\Phi$ 
being  nondegenerate, it
provides  $L(\triv)$ with
a  Frobenius  algebra structure
induced by a linear function on $\C[\h]/I$
that factors through the top degree homogeneous component. Thus,
 the scheme
$\Spec L(\triv)$ is Gorenstein (of dimension zero).~\end{proof}

\subsection{Structure of the paper}
The  paper is organized as follows. 
In \S2, we introduce an important tool
of this paper -- the BGG resolutions of irreducible 
finite dimensional representations. 
In \S3, using the Knizhnik-Zamolodchikov functor
\footnote{This functor
was originally introduced by Opdam, 
and exploited in [BEG].
For a  more detailed study of the  Knizhnik-Zamolodchikov functor 
the reader is referred to [GGOR].}
 and 
the representation theory of Iwahori-Hecke algebras
at roots of unity [DJ1,DJ2], 
we prove the existence of the BGG resolution
in the type ${\mathbf{A}}$ case. This will play a crucial step in the proof 
of  Theorem
\ref{classif}, Theorem \ref{perv}, and Character Theorem
\ref{char} (for type ${\mathbf{A}}$) given
in \S5. Section 4 contains some general
results about category $\oo(\hh_c)$,
and also  a proof of Theorem \ref{BD1}(i)  and 
Theorem \ref{char_ehe}(i)-(ii).   Classification Theorem
\ref{classif} is proved in \S5, using 
 the trace techniques developed in \cite[\S5]{BEG}.
In \S6 we prove the results on finite dimensional representations 
of rational Cherednik algebras for types ${\mathbf{B}}$ and ${\mathbf{D}}$.
This section contains a proof of Theorem \ref{BD1}(ii)-(iii).
Finally, in \S7 we discuss connections with other topics, such as 
the affine flag variety, and 
the geometry of the Hilbert scheme of the affine plane. 
\medskip

\begin{remark} After the bulk of this paper had been written, we
received a preliminary version of a very interesting
preprint by I. Gordon [Go], containing results
essentially equivalent to our Theorems \ref{BD1}(i), and 
\ref{char}, \ref{char_ehe}, and \ref{reso} (for the case $c=1+\frac{1}{h}$).
The main ideas of [Go] are very similar to ours. However, reading [Go] has allowed us 
to improve a number of the results of the present paper. 
\end{remark}

\noindent
{\bf Acknowledgements.} {\footnotesize{ We are very grateful to
 I. Cherednik for several useful discussions, and to I. Gordon
who kindly provided us with a preliminary version of [Go]
before it was made public. 
The first author was partially supported by the NSF grant 
DMS 00-71792 and A.~P.~Sloan Research Fellowship; the
work of the second author was  partly conducted for 
the Clay Mathematics Institute and partially
supported by the NSF grant DMS-9988796}.}

\section{The BGG resolution}
\setcounter{equation}{0}
\subsection{Special case: the Coxeter number}\label{sp}

{\it One dimensional} representations of the algebra
$\hh_c(W)$ are easily described by the following 
well known elementary result.

\begin{proposition}\label{coxeter}
\vi
The assignment 
$$w\mapsto \id\,,\, x\mapsto 0\,,\, y\mapsto 0\quad,\quad
(\forall w\in W, x\in \h^*, y\in \h)\,$$
 can be extended to
a 1-dimensional $\hh_c$-module if and only if the equation
$2\sum_{s\in S} c_s=\ell$ holds.

\vii If $c$ is  a {\sl constant} function, then the above
assignment can be extended to
a 1-dimensional $\hh_c$-module if and only if
$c=1/h$.
\end{proposition}
\begin{proof} The assigment above
sends the last commutation relation in \eqref{defrel} to
\begin{align*}
0\,&= y\cd x-x\cd y \,=\, \langle y,x\rangle\cd \id
-\sum\nolimits_{_{s\in S}}
c_s\langle y,\alpha_s\rangle
\langle\alpha_s^\vee,x\rangle \cd \id. 
\end{align*}
But for any conjugacy class $C$ of reflections, one has
$$\langle y,x\rangle =
\frac{\ell}{2|C|}\sum\nolimits_{_{s\in C}}\,\langle y,\alpha_s\rangle
\langle\alpha_s^\vee,x\rangle\quad,\quad\forall y\in\h\,,\,x\in\h^*
$$
(the two sides are clearly proportional, and the proportionality 
coefficient is found by substituting base elements $x=x_i$, $y=y_i$, and
summing over $i=1,\ldots,\ell$). 
Thus we see that the assigment of the Proposition gives rise to
a 1-dimensional $\hh_c$-module if and only if
we have $2\sum_{s\in S}\, c_s=\ell$.

Part (ii) is now immediate from the identity $2|S|= h\cd \ell$,
see \cite[3.18]{hum} (in the Weyl group case it goes back to Kostant, [Ko]).
\end{proof}

\begin{remark} More generally, let $\hh_c=\hh_c(V,\G,\omega)$
be a symplectic reflection algebra associated to a symplectic vector
space $(V,\omega)$ and a finite group
$\G\sset Sp(V,\omega)$ , see [EG]. Let $S\sset \G$ denote the
set of symplectic reflections. Then, a similar computation shows that
trivial 1-dimensional $\hh_c$-modules exist if and only if 
the following  linear equation on $c$ holds: $\sum_{s\in S}\,
c_s\cdot\omega_s=\omega$ (see [EG] for the notation  $\omega_s$).
$\lozenge$\end{remark}

\subsection{BGG resolutions of finite dimensional modules}
Since the 1-dimensional $\hh_c$-module considered
 in Proposition \ref{coxeter} restricts to the trivial
representation $\triv$ of the group $W$, it is isomorphic
to $L(\triv)$. We have the following Cherednik algebra analogue of the 
Bernstein-Gelfand-Gelfand resolution of the trivial
representation of a semisimple Lie algebra (it also appears in [Go]).

\begin{proposition}\label{BGG_cox} Let $\hh_c$ be the rational Cherednik
algebra associated to $W$, and let $2\sum_{s\in S} c_s=\ell$.
 Then there exists a resolution 
of  $L(\triv)$ by standard modules, of the form
$$
0\to M(\wedge^{\ell}\h^*)\to\ldots\to M(\wedge^2\h^*)\to
M(\h^*)\to M(\triv)\to L(\triv)\to 0\,.
$$
Moreover, when restricted to the polynomial subalgebra $\C[\h]\subset
\hh_c$, 
the resolution above becomes the standard Koszul complex:
$$0\to \C[\h]\otimes\wedge^{\ell}\h^*\to\ldots\to
\C[\h]\otimes\wedge^2\h^*\to\C[\h]\otimes\h^*\to\C[\h]\to\C\to 0\,.
$$ 
\end{proposition}
\begin{proof}
We build the required resolution 
explicitly by induction.
The first nontrivial arrow is the augmentation $ M(\triv)=\C[\h]\onto
L(\triv)= \C$.
Its kernel is the space of 
polynomials vanishing at $0$, which is an $\hh_c$-submodule,
generated (already over $\C[\h]$) by linear polynomials, which 
are singular vectors in this module.
Thus, we have $\hh_c$-module morphisms:
$M(\h^*)\to M(\triv)\to L(\triv)=\C$.
 Now, we claim that the kernel of the first map
is a submodule of $M(\h^*)$, which is generated (already over $\C[\h]$)
by $\wedge^2\h^*$ sitting in degree 1 of $M(\h^*)$. 
This is clear because, from commutative algebra
point of view, we have just the first two terms of the Koszul
resolution of $\C$.
Thus, we obtain  $\hh_c$-module morphisms: $
M(\wedge^2\h^*)\to M(\h^*)\to M(\triv)\to L(\triv)$.
Continuing in the same fashion, we get the desired resolution,
which in this case coincides with the Koszul resolution of $\C$ as a 
$\C[\h]$-module.
\end{proof}

In the special case of the root system of type $\mathbf{A_{n-1}}$, 
the analogy with the BGG-resolution extends  further, to {\it all}
finite dimensional $\hh_c(S_n)$-modules:

\begin{theorem}[BGG-resolution]\label{reso}
Assume $R=\mathbf{A_{n-1}}$. If  $c=r/n>0$, where $(r,n)=1\,,\,r\in \Z,$ 
then $L(\triv)$ admits a resolution 
 of the form
$$
0\to M(\wedge^{\ell}\h^*)\to\ldots\to M(\wedge^2\h^*)\to
M(\h^*)\to M(\triv)\to L(\triv)\to 0\,.
$$
\end{theorem}

The proof of Theorem \ref{reso} is given in Section 3. 
Theorem \ref{reso} will be used in the proof of Theorem \ref{classif}.

Similarly, outside of type $\mathbf A$, we have 
the following (weaker) generalization of Proposition
\ref{BGG_cox}.

\begin{theorem}[BGG-resolution]\label{reso1}
Let $W$ be a finite Weyl group and `$c$' be
 a constant of the form $\frac{1}{h}+k$, where $k=0,1,2,\ldots$.
 Then $L(\triv)$ admits a resolution 
 of the form
$$
0\to M(\wedge^{\ell}\h^*)\to\ldots\to M(\wedge^2\h^*)\to
M(\h^*)\to M(\triv)\to L(\triv)\to 0\,.
$$
\end{theorem}

For $k=1$ this theorem is contained in [Go]. 
The proof of Theorem \ref{reso1} is given in Section 4.
It will be used for the proof of Theorem \ref{BD1}(i).

\section{The Hecke algebra
and the proof of Theorem \ref{reso}}
\setcounter{equation}{0}
\subsection{The Hecke algebra and Specht modules}
Recall that to a finite Coxeter group $(W, S)$,
and  a $W$-invariant function $q:S\to \C^*$ 
one associates
the corresponding
 Hecke algebra $\mathcal H_q=\mathcal H_q(W)$.
In this subsection we will recall from \cite{DJ1,DJ2}
some known facts 
about the representation theory of $\mathcal H_q$
for root system of type $\mathbf{A}$
 which will be used below. 

\begin{remark} We adopt the normalisation in which the
quadratic relations for the generators
$\{T_i=T_{s_i}\}_{i=1,\ldots,\ell}$
of our  Hecke algebra $\mathcal H_q=\mathcal H_q(W)$
read: $(T_i-1)(T_i+q)=0$.
 In \cite{DJ1,DJ2}, the quadratic relations for $T_i$ are 
$(T_i+1)(T_i-q)=0$, so the $T_i$ of \cite{DJ1} are different from ours 
by a factor $q$, and our $q$ corresponds to $q^{-1}$ of 
\cite{DJ1,DJ2}. 
$\lozenge$\end{remark}

Let $e$ be the order of $q$ (equal to infinity if $q$ is not a root of $1$).   
Recall [DJ1] that for every partition $\lambda$ of $n$, 
we have the Specht module $S^\lambda=S^\lambda(q)\subset \mathcal H_q$ 
over $\mathcal H_q$, 
which is flat with respect to $q$, viewed as a parameter. We also have  
its quotient $D^\lambda$,
which is nonzero  if and only if  $\lambda$ is $e$-regular (i.e. multiplicities of 
all parts of the partition $\lambda$
are smaller than $e$), and is irreducible when nonzero; 
this gives the full list of irreducible $\mathcal H_q$-modules without 
repetitions. Moreover, by \cite[7.6]{DJ1}, the module $D^\lambda$ occurs in $S^\lambda$ 
with multiplicity $1$. 

Recall also from [DJ2] that $S^\lambda$ and $S^\mu$ are in the same block 
if and only if $\lambda$ and $\mu$ have the same $e$-core 
(see \cite{JK} for definition of core). In particular, if 
$\lambda$ is an $e$-core itself (i.e., has no hook of length $e$; 
in this case it is also $e$-regular),
then $S^\lambda$ lies in a block separate from other $S^\mu$ and, 
furthermore, 
$D^\lambda\ne 0$. Hence, if $\lambda$ is an $e$-core, then
the only irreducible module in the block of 
$S^\lambda$ is $D^\lambda$, and consequently $S^\lambda$ is simple. 

\subsection{The Hecke algebra at the primitive n-th root of 1}
The simplest nontrivial special case $e=n$ 
(i.e., $q$ is a primitive $n$-th root of unity) will be especially important
for us. In this case, $n$-core diagrams are all but single-hook
diagrams, i.e. diagrams 
$\lambda^i$ corresponding 
to representations $\wedge^i\h$ of $S_n$, $i=0,...,n-1$. 
The $n$-core of $\lambda^i$ is the empty diagram. 
Thus all Specht modules $S^{\lambda^i}$ belong to the same block, 
while other Specht modules are all irreducible. 

Moreover, it is seen from \cite{DJ1} that if $\lambda\ne \lambda^i,
\forall i,$ then 
$S^\lambda$ is projective.
Therefore, the block of $S^\lambda$ consists of multiples of $S^\lambda$. 
Thus, we have: 
$$
\mathcal H_q=\bigl(\oplus_{\lambda \text{ is an $e$-core}}\;{\rm End}_{_\C} S^\lambda\bigr)
\oplus A,
$$ 
 where $A$ is an indecomposable, non-semisimple algebra, 
corresponding to the only non-semisimple block of the representation category
of $\mathcal H_q$. 

Further, we need a more precise description of the Specht modules
for $\lambda=\lambda^i$, which can be inferred from \cite{DJ1,DJ2}. 
First of all, $D^{\lambda^i}$ is nonzero  if and only if  $i\ne n-1$. 

\begin{proposition}\label{exactseq}
The Specht module $S^{\lambda^i}$ 
has composition factors $D^{\lambda^i}$, $D^{\lambda^{i-1}}$
(only one of them if the other is zero or not defined). 
\end{proposition}
\noindent
(Here and elsewhere, each simple composition factor is written as many
times
as it occurs in the composition series.)
\subsection{Knizhnik-Zamolodchikov functor on standard modules}
Abusing notation, we will write $M(\lambda)$ for 
the standard module over $\hh_c$ whose lowest weight is the representstaion 
of $S_n$ corresponding to a partition $\lambda$. 

Let $c$ be any complex number, and 
 $q=e^{2\pi ic}$.
Recall (\cite{OR}, see also \cite{BEG})
 that for any $V\in \oo(\hh_c)$,
its localization to $\hreg:= \h\ssminus (\cup_{s\in S}\,H_s)$
(= set of regular elements), is a vector bundle $V\mr$ with flat
connection,
i.e., a local system. Assigning to $V$ the monodromy of 
that local system gives  the ``Knizhnik-Zamolodchikov functor'' 
$\kz: \oo(\hh_c)\to {\rm Rep}(\mathcal H_q)$, see   [GGOR]
for more details.
Put $\tilde S^\lambda:=\kz(M(\lambda))$.

\begin{lemma}\label{Specht}
The $\mathcal H_q$-modules 
$\tilde S^\lambda$ and $S^\lambda$ have the same composition factors.
\end{lemma}

\begin{proof}
We begin with the following standard result (see e.g. \cite[Lemma 2.3.4]{CG}).
\smallskip

\noindent
{\bf Claim:\;} {\it
Let $A$ 
be a flat, finite rank algebra over $\C[[t]]$.
Let $M,N$ be two $\C[[t]]$-flat, finite rank $A$-modules. 
Suppose that $M\otimes_{\C[[t]]}\C((t))$ and $N\otimes_{\C[[t]]}\C((t))$ 
are isomorphic as $A\otimes_{\C[[t]]}\C((t))$-modules. 
Then $M/tM$ and $N/tN$ have the same composition factors as 
$A/tA$-modules.}  \qed
\smallskip

Now, in our particular situation, we have two flat (holomorphic) families 
of modules over $\mathcal H_q$ (where $q=e^{2\pi i c}$),
namely $S^\lambda(c)$ and $\tilde S^\lambda(c)$.
They are obviously 
isomorphic at $c=0$, and hence for regular $c$, by a standard deformation argument. 
Now, the Claim above implies
that, for special value of $c$, the modules
$S^\lambda$ and $\tilde S^\lambda$ have the same composition factors. 
Lemma \ref{Specht} is proved.
\end{proof}

\begin{corollary} \label{Specht1} For $c=r/n$, where $(r,n)=1\,,\,r\in \Z,$ we have:

1) If $\lambda\ne \lambda^i\,,\,\forall i$, then $\tilde S^\lambda$ is isomorphic to 
the Specht module $S^\lambda$, and in particular is irreducible. 

2) If $\lambda=\lambda^i$, then $\tilde S^\lambda$ belongs 
to the non-semisimple block and has composition factors 
$D^{\lambda^i}$ and $D^{\lambda^{i-1}}$ (only one of them if the other
is zero). 
\end{corollary}

\begin{proof} This follows from Lemma \ref{Specht} and 
Proposition \ref{exactseq}, since $q=e^{2\pi ic}$ is a primitive 
$n$-th root of unity.   
\end{proof} 

\subsection{Construction of the BGG resolution}

Throughout this subsection let $c=r/n$, where $(r,n)=1\,,\,
r\in\Bbb N$. 

\begin{proposition}\label{structO}
\vi If $\lambda\ne \lambda^i\,,\,\forall i,$ then $M(\lambda)$ is irreducible.

\noindent
\vii If $\lambda\ne \mu$ are Young diagrams, and 
${\rm Hom}(M(\lambda),M(\mu))\ne 0$. 
Then $\lambda=\lambda^i$ and~$\mu=\lambda^{i-1}$. 
\end{proposition}

\begin{proof} According to
a result by Opdam-Rouquier presented in (\cite[Lemma 2.10]{BEG}),
the Knizhnik-Zamolodchikov functor is faithful when restricted to 
standard modules. Therefore, statement (i) follows from 
the fact, if
 $\lambda\ne \lambda^i\,,\,\forall i,$ then
 $\tilde S^\lambda$ does not belong to the block of
$\tilde S^\mu$, for $\mu\ne \lambda$, which is a consequence of 
Corollary \ref{Specht1}. 

To prove statement (ii), 
we observe that again by Corollary \ref{Specht1}, 
$\lambda=\lambda^i$ and $\mu=\lambda^j$ 
for some $i$ and $j$. 
If $j>i$, the statement is clear since one knows, see e.g. formula
\eqref{integer_i}  of \S4 below, that:
$\kappa(c,\wedge^j\h)>\kappa(c,\wedge^i\h)$.
On the other hand, if $j<i-1$, then any nonzero homomorphism 
from $M(\lambda^i)$ to $M(\lambda^j)$ would give rise, through the Knizhnik-Zamolodchikov
functor and \cite[Lemma 2.10]{BEG}, to a nonzero homomorphism from 
$\tilde S^{\lambda^i}$ to $\tilde S^{\lambda^j}$. 
But such a homomorphism does not exist, since 
these modules have no composition factors in common. 
\end{proof}

\begin{proposition}\label{Homs}
The space ${\rm Hom}(M({\lambda^i}),M({\lambda^{i-1}}))$
is one-dimensional.
\end{proposition}

\begin{proof}
The modules $\tilde S^{\lambda^i}$ and $\tilde S^{\lambda^{i-1}}$ 
have one composition factor in common, so we see that 
${\rm Hom}(M({\lambda^i}),M({\lambda^{i-1}}))$
is at most one dimensional. Thus we are left 
to show that this space is nonzero.  

By Deligne's criterion, the connection defined by Dunkl operators 
has regular singularities. Therefore, by the Riemann-Hilbert correspondence, 
the $(\hh_c)\mr$-module $M(\lambda^j)\mr$
(where ``$(-)\mr$'' denotes localization to the Zariski open
affine set 
of regular points in $\h$) is reducible for $j=1,..,n-2$, 
since so is, $\tilde S^{\lambda^j}$,  its image under $\kz$.
Hence, $M(\lambda^j)$ itself is reducible
(it has a nonzero proper submodule $M(\lambda^j)\cap N$, where $N$ 
is a nonzero proper submodule of $M(\lambda^j)\mr$). 
This implies that another standard module 
maps to $M(\lambda^j)$ nontrivially. This yields the result for $i>1$. 

It remains to prove the statement for $i=1$
(and in particular to show that $M(\triv)$ is reducible). To do this, observe that the 
validity of the statement for $i>1$ implies 
that, for any $i\ge 1$, the module  $\tilde S^{\lambda^i}$
 is included in an exact sequence 
$$
0\to D_{\lambda^i}\to \tilde S_{\lambda^i}\to D_{\lambda^{i-1}}\to 0
$$
(to afford maps $\tilde S^{\lambda^i}\to \tilde S^{\lambda^{i-1}}$
coming from $M(\lambda^i)\to M(\lambda^{i-1})$, $i>1$). 
Hence, we have a surjective map $\tilde S^{\lambda^1}\onto \tilde 
S^{\lambda^0}$. 
By the Riemann-Hilbert correspondence, 
this map gives rise to a nonzero map 
$f: M(\lambda^1)\mr\to M(\lambda^0)\mr$. 
Let $N\subset M(\lambda^1)\mr$ be defined by 
$N=f^{-1}(M(\lambda^0))$. This is a nonzero $\hh_c$-module. 
Let $P=N\cap M(\lambda^1)$. Then $P$ is a nonzero submodule of 
$M(\lambda^1)$ which maps nontrivially to $M(\lambda^0)$. 
This means that for some $\tau\in\irrep(S_n)$ we have a diagram
of nonzero maps of $\hh_c$-modules:
$M(\lambda^1)\leftarrow M(\tau)\rightarrow M(\lambda^0)$. 
But it follows from the above that the only choice for 
$\tau$ is $\tau=\lambda^1$. Thus, we have a nonzero map 
$M(\lambda^1)\to M(\lambda^0)$. We are done. 
\end{proof}

Propositions \ref{structO}
and \ref{Homs} imply the following corollary (we will often identify
$\h^*\simeq\h$).

\begin{corollary}\label{comp}
For $c=r/n>0$, $(r,n)=1\,,\,r\in{\mathbb{N}}$, there exists a (unique up to scaling) 
complex of $\hh_c$-modules 
\begin{equation}\label{eqcomp}
0\too M(\wedge^{n-1}\h)\too \ldots\too M(\h)\too M(\triv)\too
 L(\triv)\too 0.
\end{equation}
\end{corollary} 

\subsection{Exactness of the BGG resolution}
The following proposition implies Theorem \ref{reso}.

\begin{proposition}\label{comp1}
The complex (\ref{eqcomp}) is exact. 
\end{proposition}

To prove the Proposition, we introduce the notion
of  {\it thick} object of $\oo(\hh_c)$. 
We  say that $V\in \oo(\hh_c)$ is thick if 
$V\mr\ne 0$. Otherwise we say that $V$ is {\it thin}. 
In other words, $V$ is thick  if and only if  its Gelfand-Kirillov
dimension 
equals $n-1$, and thin  if and only if  it is $< n-1$. 
 
\begin{lemma}\label{thinthick}
If $c=r/n>0$, $(r,n)=1$, then 
$L(\tau)$ is thin if and only if $\tau=\triv$. 
\end{lemma}

This lemma
implies in particular that all $L(\tau)$ except $L(\triv)$ are infinite 
dimensional. 

\begin{proof}
First of all, $L(\triv)$ is thin. Indeed, we have an exact sequence 
$M(\lambda^1)\to M(\lambda^0)\to V\to 0$, 
and $L(\triv)$ is a quotient of $V$. 
Applying to this sequence the Knizhnik-Zamolodchikov functor 
and using the fact that this functor is exact, 
we get an exact sequence
$\tilde S^{\lambda^1}\to \tilde S^{\lambda^0}\to \kz(V)\to 0$.
Thus, $\kz(V)=0$, so $V$ and hence $L(\triv)$ is thin.

Now, we claim that all other $L(\tau)$ are thick. Indeed, by the Riemann-Hilbert 
correspondence, 
$\kz(L(\tau))$ is either zero or irreducible, and 
the set of all irreducibles obtained this way 
has to be the full set of irreducibles $D^\lambda$ of $\mathcal H_q$.
But the number of irreducibles in our case is $p(n)-1$, where $p(n)$ is 
the number of partitions of $n$. Hence, for all  $\tau\ne \triv$, the module
$\kz(L(\tau))$ has to be nonzero.
Thus, $L(\tau)$ is thick, for any  $\tau\ne \triv$.
\end{proof}

\begin{corollary}\label{length2}
For any $i=0,1,...,n-2$, the module $M(\lambda^i)$ 
is included in an exact sequence 
$$
0\to L(\lambda^{i+1})\to M(\lambda^i)\to L(\lambda^i)\to 0. 
$$
\end{corollary}

\begin{proof} Let  $J(\lambda^i)$ be the maximal proper submodule 
of $M(\lambda^i)$. By  exactness of the functor
$\kz$ we get $\kz(J(\lambda^i))=D^{\lambda^{i+1}}$.
Since any $L(\tau)$ for $\tau\ne \triv$ is
 thick, by Lemma \ref{thinthick}, we see that  $J(\lambda^i)$
must have a  composition factor $L(\lambda^{i+1})$
with possibly other 
composition factors all isomorphic to  $L(\triv)$.
But  $L(\triv)$ cannot occur in   $J(\lambda^i)$
since $\kappa(c,\triv)$, see \eqref{kapa},
 is smaller than the lowest 
eigenvalue of $\bold h$ in $J(\lambda^i)$. Thus,
$J(\lambda^i)\simeq L(\lambda^{i+1})$.
\end{proof}

Now, Proposition \ref{comp1} easily follows from 
Corollary \ref{length2}. 

\section{Results on category $\oo(\hh_c)$}
\setcounter{equation}{0}
\subsection{${\mathtt{p}}$-function and characters}
 Let $(W,S)$ be a Coxeter
group, $C\subset S$ a conjugacy class of reflections,
and  $\tau\in \irrep(W)$. For $s\in S$, the equation $s^2=1$ 
implies  $\tau(s)^2=\id$, hence the linear map
 $\tau(s)$ has eigenvalues $\pm 1$. For any 
 $s\in C$, let $q$ be the number of
$(+1)$-eigenvalues, and $p$ the number of
$(-1)$-eigenvalues  in  $\tau(s)$. Then, $q+p=\dim\tau$, and also
$\Tr(s, \tau)=q-p$. Now,  the central element $\sum_{s\in C} \,s$ acts
in  $\tau$ as a scalar, say ${\mathtt{p}}(C,\tau)$.
We compute: 
$${\mathtt{p}}(C,\tau)\cdot\dim\tau=
\Tr\bigl(\sum\nolimits_{s\in C}
\,s,\tau\bigr)=|C|\cdot\Tr(s, \tau)=|C|\cdot(q-p)=|C|\cdot(\dim\tau-2p)\,.
$$
On the other hand, the scalar  ${\mathtt{p}}(C,\tau)$
is well-known to be an algebraic integer, for  any conjugacy
class $C$
in a finite group
and for any irreducible module $\tau$. Thus,
we conclude 
\begin{equation}\label{integer}
{\mathtt{p}}(C,\tau)\;=\;\frac{|C|\cdot(q-p)}{\dim\tau}\;=\;
|C|-\frac{2|C|}{\dim\tau}\cdot p\quad\text{{\sl is an} {\sf integer}}\,.
\end{equation}

Write ${\mathtt{p}}(\tau):=\sum_{C\subset
S}\,{\mathtt{p}}(C,\tau)$ 
for the  scalar by which
  the central element $\sum_{s\in S} \,s$ acts
in  $\tau$. From \eqref{integer} we deduce that ${\mathtt{p}}(\tau)$ is
an integer, moreover, we have:
\begin{equation}\label{inequality}
-|S|={\mathtt{p}}(\sign)\, <\, {\mathtt{p}}(\tau)\,<\, {\mathtt{p}}(\triv)=|S|
\quad,\quad\forall\tau\in\irrep(W)\,,\,\tau\neq\triv,\sign.
\end{equation}
In the special
case $\tau=\wedge^i\h$ that will play a role below,
  one easily calculates using  \eqref{integer}:
\begin{equation}\label{integer_i}
{\mathtt{p}}(\wedge^i\h) =|S|-h\cdot i \,,
\quad\text{hence}\quad \kappa(c,\wedge^i\h)= \ell/2 -c\cdot{\mathtt{p}}(\wedge^i\h)=
\ell/2+c\cd h(i - \ell/2).
\end{equation}
(This follows also from Proposition \ref{BGG_cox} 
because the degrees in the consecutive terms of the
 Koszul resolution go with step 1.)

Next, let  $K(\oo(\hh_c))$ be  the Grothendieck group of the
category  $\oo(\hh_c))$.
Recall that classes
$\big\{[M(\tau)]\in K(\oo(\hh_c))\,,\,\tau\in \irrep(W)\big\}$ form a $\Z$-basis
in $K(\oo(\hh_c))$. 
Thus, for any  object $V$ there are certain (uniquely
determined) integers $a_\tau\in\Z$ such that in 
 $K(\oo(\hh_c))$ we have: $\dis [V]=\sum\nolimits_{\tau\in
\irrep(W)}\,a_\tau\cdot [M(\tau)].$ Now, assuming
 `$c$' is a constant function on $S$, from formula \eqref{kapa}
for the character of standard modules we get

\begin{equation}\label{char_L}
\Tr_{_{V}}(w\cdot  t^{\bold h})
=\frac{t^{\ell/2}}{\det(1-t\cdot w)}\cdot\sum\nolimits_{\tau\in
\irrep(W)}\,a_\tau\cdot\Tr(w,\tau)\cdot t^{-c\cdot{\mathtt{p}}(\tau)}\,.
\end{equation}
In particular, we find
\begin{equation}\label{char_L2}
\Tr_{_{V}}(t^{\bold h})
=\frac{t^{\ell/2}}{(1-t)^\ell}\cdot\sum\nolimits_{\tau\in
\irrep(W)}\,a_\tau\cdot\dim\tau\cdot t^{-c\cdot{\mathtt{p}}(\tau)}\,.
\end{equation}

\subsection{Proof of the character  Theorem \ref{char} and  Theorem
\ref{BD1}(uniqueness)}
Without loss of generality, we may (and will)
restrict our attention to the case $\varepsilon=1$,
i.e. $c>0$. The assumptions of  Theorem
\ref{BD1} imply that $(r,h)=1$ in all 
the cases at hand. Hence,
$q=e^{2\pi i c}$ is a primitive root of unity of
order $h$. Hence, by \cite[Lemma 4.3]{Go}, $\oo(\hh_c)=\oo'\oplus
\oo'',$ where $\oo'$ is generated
by the $L(\wedge^i\h)\,,\,i=0,\ldots,\ell-1,$ and $\oo''$ is generated
by the other simples in $\oo(\hh_c)$. Thus, the simple and standard modules
in
$\oo''$ coincide, hence no simple module in $\oo''$ can be finite
dimensional.
 
Now, fix an irreducible  finite dimensional $\hh_c$-module $L$.
\smallskip

\noindent
{\sl Proof of the character  Theorem \ref{char}.\;} Our argument
is similar in spirit to H. Weyl's proof of his character
formula for a compact Lie group.

First of all, as we have seen above,  $L$ must be an object of  $\oo'$.
Therefore in the Grothendieck group of the
category  $\oo'$ we can write
$[L]=\sum a_i\cdot[M(\wedge^i\h)],$ for some integers $a_i$.
Hence, applying a version of formula \eqref{char_L2} in our specific
situation
and using \eqref{integer_i},
we find:  $\Tr_L(t^{\mathbf h})=
t^{-\ell(c\cdot h-1)/2}\cdot Q(t^{c\cdot h})/(1-t)^\ell,$ where
$Q(z)=\sum_{i=0}^\ell\, z^i\cdot \left({\ell\atop i}\right)\cdot a_i$.

The polynomial
$Q\in \C[z]$ has degree $\ell,$ and the ratio $Q(t^{c\cdot
h})/(1-t)^\ell$
cannot have a pole
at $t=1$ since it is equal to $\dim L$, which is finite. We deduce
that $Q(z)=m\cdot (1-z)^\ell,$ for some 
constant $m\in\C$. It follows  that, for any $i$ we have: $a_i=m\cdot(-1)^i$.
Thus we find: $[L]= m\cdot\sum_{i=0}^\ell\,(-1)^i\cdot [M(\wedge^i\h)]$.
This implies, since $a_i\in\Z,$ that the constant $m$ must be an integer.
But if the class  $[L]\in K(\oo')$ of a nonzero object  $L\in\oo'$
is divisible by an integer $m$, then
 the object $L$ cannot be simple unless $m=\pm 1$.
Hence the irreducibility of $L$ forces $m=\pm 1$.
Therefore, for the character of $L$ we get $\Tr_L(t^{\mathbf h})=\pm
t^{-\ell(c\cdot h-1)/2}\cdot(1-t^r)^\ell/(1-t)^\ell$. Letting $t\to 1$ 
here we must get the dimension of $L$, which is positive.
Therefore, one must have
the `$+$'-sign in the formula.
Thus,  $[L]=\sum_{i=0}^\ell\,(-1)^i\cdot [M(\wedge^i\h)]$. This proves
that the character of $L$ is given by the formula of Theorem \ref{char}.

It remains to show that $L\simeq L(\triv)$. Suppose
$L\simeq L(\wedge^j\h),$ for some $j>0$.
Then, since $\kappa(\wedge^j\h) > \kappa(\triv)$ by \eqref{integer_i}, the power
expansion: $\Tr_L(t^{\mathbf h})=\sum\, (\dim L_a)\cdot t^a$ 
cannot contain $t^{\kappa(\triv)}$ with a  non-zero coefficient.
But on the other hand, it is clear e.g. from formula
\eqref{char_L2} that  $t^{\kappa(\triv)}$
does appear in  the power
expansion of $\sum_{i=0}^\ell (-1)^i\cdot \Tr_{_{M(\wedge^i)}}(t^{\mathbf
h})$
with  a  non-zero coefficient. The contradiction forces $j=0$, and 
 $L\simeq L(\triv)$. 

This completes the proof of 
Theorem  \ref{char} and at the same time yields the uniqueness claim
in Theorem
\ref{BD1}.\qed

\subsection{Morita functors}\label{shift}
Let $\varepsilon: W \to \C^\times$ be a character of $W$, and
$\bold 1_\varepsilon$  the characteristic function of the
subset $\{s\in S\;|\;\varepsilon(s)=-1\}$.
We write $\e_\varepsilon=\frac{1}{|W|}\sum_{w\in W}\,\varepsilon(w)\cdot
w$ for
 the central idempotent in $\Bbb C{W}$ corresponding to $\varepsilon$. 

\begin{proposition}\label{isom} There exists a filtered algebra
isomorphism
$\phi_c: \e \hh_c\e\iso \e_\varepsilon\hh_{c+\bold 1_\varepsilon}\e_\varepsilon$. 
\end{proposition}

\begin{proof} Let $A_1(c)$ denote the first algebra, 
and $A_2(c)$ the second one. 
Let $\dd(\hreg)$ 
denote the algebra of differential operators on 
the regular part of $\h$. 
The algebras $A_1(c)\,,\, A_2(c),$ each carry a filtration 
given by $\deg(\h^*)=\deg(w)=0$, $\deg(\h)=1$, and a grading 
given by $\deg(\h^*)=1$, $\deg(\h)=-1$, $\deg(w)=0$. 
Similarly, $\dd(\hreg)$ has a filtration 
by order  of differential operators,
and grading by their homogeneity degree. 

As explained in [EG], [BEG]
we have the spherical homomorphism
$\theta_c: A_1(c)\to \dd(\hreg)$, and
`$\varepsilon$-antispherical' homomorphism
$\theta_c^-: A_2(c)  \to  \dd(\hreg)$. 

It is shown in [BEG]
that for generic $c$ (i.e. all but countably many), 
the images of these two homomorphisms are the
same. Namely,  for an
orthonormal basis $y_i$ of $\h$, the elements $\e(\sum y_i^2)\e\in A_1(c)$ and
$\e_\varepsilon(\sum y_i^2)
\e_\varepsilon\in A_2(c)$,
go to the same operator, the so-called Calogero-Moser operator
${\mathbf{L}}_c$
(see [EG]);
the same holds for the elements $\e P(x)\e$ and ${\emi}P(x){\emi}$, where 
$P$ is a polynomial on $\h$.
On the other hand, 
$A_1(c)$ is generated by $\e(\sum y_i^2)\e$ and $\e\cdot\C[\h]\cdot\e$, and
similarly for 
$A_2(c)$. This proves:
${\mathtt{Image}}(\theta_c)={\mathtt{Image}}(\theta_c^-)$,
for generic `$c$'. We are going to show that the equality holds, in
effect, for
{\it any}  `$c$'.

To this end, 
recall that the parameter  `$c$' varies over an affine space $\C^d$.
It is known, see [EG], that for any $c\in \C^d$,
both $\theta_c$ and $\theta_c^-$ are injective filtration
preserving morphisms, such that the associated graded maps are also
injective. Moreover, the dimensions of associated graded components
of ${\mathtt{Image}}(\theta_c)$, resp. of ${\mathtt{Image}}(\theta_c^-)$,
are independent of  `$c$'. 
Thus,
as `$c$' runs over  $\C^d$,
 one may treat the spaces  ${\mathtt{Image}}(\theta_c)$, resp.
 ${\mathtt{Image}}(\theta_c^-)$, as the fibers of a
filtered subbundle ${\mathbf{Image}}(\theta^+)$, resp.
 ${\mathbf{Image}}(\theta^-)$, in a
{\it trivial} vector bundle on  $\C^d$ with fiber $\dd(\hreg),$
an infinite dimensional filtered vector space.

 By an argument two paragraphs above we know that, on a dense subset
$U\subset \C^d$, the two vector
subbundles coincide: ${\mathbf{Image}}(\theta^+)|_U=
{\mathbf{Image}}(\theta^-)|_U.$
It follows that they coincide everywhere, i.e., 
${\mathtt{Image}}(\theta_c)={\mathtt{Image}}(\theta_c^-)$,
for any $c\in \C^d$. Thus the  filtered isomorphism claimed by the 
Proposition may be given by
$\phi_c:=(\theta_c^-)^{-1}\ccirc\theta_c$.
\end{proof}

Let $C$ be the set of those $c$ 
for which $\hh_c\e V=V$ for any 
$V\in \mathcal \oo(\hh_c)$. 
This condition is equivalent to saying that 
the functors $F_c:\oo(\hh_c)\to \oo(\e \hh_c\e)$ and 
$G_c: \oo(\e \hh_c\e)\to \oo(\hh_c)$ 
given by $F_c(V)= \e V$, 
$G_c(Y)=\hh_c\e\otimes _{\e \hh_c\e}Y$ 
are mutually inverse equivalences of categories. 
It is easy to see that $c\in C$  if and only if  
$\hh_c\cdot\e\cdot V=V$ for all irreducible objects of $\oo(\hh_c)$.

The isomorphism $\phi_c$ of Proposition \ref{isom}
induces an equivalence 
$\Phi_c: \oo(\ehe)\iso$
$
\oo({\e_\varepsilon}\hh_{c+{\mathbf 1}_\varepsilon} {\e_\varepsilon})$. 

\begin{lemma}\label{morita} 
If $c$ and $-c-{\mathbf 1}_\varepsilon$ belong to $C$ then 
the shift functor $\Bbb S_{c,\varepsilon}:\oo(\hh_c)\to
\oo(\hh_{c+{\mathbf 1}_\varepsilon}),$
 given by: 
$V \mapsto \Bbb S_{c,\varepsilon}(V)=
\hh_{c+{\mathbf 1}_\varepsilon}{\e_\varepsilon}
\otimes_{ {\e_\varepsilon} \hh_{c+{\mathbf
      1}_\varepsilon}{\e_\varepsilon} }\,
\Phi_c(\e V)$
is an equivalence of categories. 
\end{lemma}

\begin{proof} The functor $\Bbb S_{c,\varepsilon}$ is a composition of three functors:
$\Bbb S_{c,\varepsilon}=G_c^- \ccirc \Phi_c\ccirc  F_c$,
where
$G_c^-(Y):=\hh_{c+{\mathbf 1}_\varepsilon}{\e_\varepsilon}\otimes_{{\e_\varepsilon}\hh_{c+
{\mathbf 1}_\varepsilon}{\e_\varepsilon}}Y$.
The functor $F_c$ is an equivalence 
since $c\in C$, and the functor $\Phi_c$ is an equivalence 
by definition. To show that $G_c^-: \oo(\e_\varepsilon\hh_{c+{\mathbf
1}_\varepsilon}\e_\varepsilon)\to\oo(\hh_{c+{\mathbf 1}_\varepsilon})
$ is an equivalence, 
recall that there exists an isomorphism $\hh_c\simeq \hh_{-c}$ which acts trivially 
on $\h$ and $\h^*$ and maps any reflection $s\in S$ to $-s$. 
Applying this isomorphism and using that $-c-{\mathbf
  1}_\varepsilon\in C$, 
we get the result. 
\end{proof} 

\begin{proposition}\label{sufflarge} Let $c$ be a constant
and $\varepsilon$ be the sign character. 
The shift functor $\Bbb S_{c,\varepsilon}: \oo(\hh_c)\to
\oo(\hh_{c+1})$
is an equivalence of categories in either of the following 
two cases

\vi If $|c|$ is sufficiently large;

\vii (see \cite{Go})\;\; If $c=\frac{1}{h}+k$, where $k=0,1,2,\ldots$.
\end{proposition} 

\begin{proof}  To handle case (i), let
 $|c|$ be large enough. Then for any $\sigma,\tau\in \irrep(S_n)$, 
either 
$\kappa(c,\sigma)=\kappa(c,\tau)$ or $|\kappa(c,\sigma)-\kappa(c,\tau)|=
c\cdot|{\mathtt{p}}(\tau)-{\mathtt{p}}(\sigma)|$ is large. 
This shows that for any $\tau$, the $\bh$-weight 
spaces of $M(\tau)$ and $L(\tau)$ 
coincide up to a large degree. Thus, $L(\tau)$ contains an invariant, and 
hence is generated by its invariants, as desired. 

To prove (ii), 
we will use \cite[Lemma 4.3]{Go}, which says in particular that 
$M(\tau)$ is irreducible unless $\tau=\wedge^i\h$
(this Lemma is stated in [Go] for Weyl groups, but remains true for a
general Coxeter group, see \cite[Remark 6.9]{BGK}). 
By this lemma, it suffices to show that 
if $c>0$ then $L(\wedge^i\h)$ is generated by its invariants,
and furthermore if $c>1$ then $L(\wedge^i\h)$ is also
generated by its
anti-invariants. In other words, we must show that 
$\e L(\wedge^i\h)\ne 0$  
for $c>0$, and (using the isomorphism $\ehe\simeq \e\hh_{-c-1}\e,$
see [BEG2])  that $\e L(\wedge^i\h)\ne 0$ for $c<-1$.  

Using Lemma \ref{sollem} and the formula  for
$\kappa(c,\wedge^i\h)$ given in \eqref{integer_i},
we obtain: 
 $\dis\Tr_{_{\e M(\wedge^i\h)}}(t^\bh)$
$\dis=t^{(1-c\cdot h)\cdot\ell/2}\cdot
\prod\nolimits_{i=1}^\ell \,(1-t^{m_i+1})^{-1}\cdot\sum
t^{c\cdot h\cdot i+m_{j_1}+...+m_{j_i}},$
where $m_i=d_i-1$ denote the {\it exponents} of $W$.
For $c>0$, we see that the powers of $t$ in the first nonvanishing term in the formal series 
expansion of the function above are increasing with $i$.
Therefore, these functions are linearly independent. Hence the result.
For $c<-1,$ the lowest terms are decreasing (note that
$m_j<h\,,\,\forall j$). Hence,
they are also independent. (This argument is due to
I.Gordon, \cite[Lemmas 4.4, 4.5]{Go}).
\end{proof} 

\section{Classification of finite-dimensional modules in type
$\mathbf A$}
\setcounter{equation}{0}
\subsection{The trace formula}\label{trace_sub}
In this subsection, we consider the case $W=S_n$. 
Let $\Tr: \hh_c\to \hh_c/[\hh_c,\hh_c]$ be the tautological projection. 
It follows from \cite[Theorem 1.8(i)]{EG}, that the space 
$\hh_c/[\hh_c,\hh_c]$ is 1-dimensional for all $c$ except possibly
countably many. 
Also, as has been shown in \S5 of \cite{BEG}, 
for all `$c$' but at most countably many, one has $\Tr(1_{_{\hh_c}})\ne 0$.
Hence, there exists at most a countable set $T\subset \C$ 
such that one has: {\it for any fixed noncommutative 
polynomial $P$ in the generators $w\in S_n$, $\{x_i\},$ and $\{y_i\}$ (dual bases
in $\h^*$ and $\h$ respectively),
there exists a function $f_P:\C\ssminus T\to \C$ 
such that}:
$\Tr(P)=f_P(c)\cdot\Tr(1_{_{\hh_c}})\;(\textsl{mod}\,
[\hh_c,\hh_c])$,
$\forall c\in \C\ssminus T$. 

\begin{lemma}\label{rati}
The function $f_P$ is rational. 
Moreover, for all 
but finitely many values of $c$,
one has
 $\Tr(P-f_P(c)\cdot 1_{_{\hh_c}})=0$.
\end{lemma}

\begin{proof} Let $F_\bullet \hh_c$ be the standard increasing
filtration of $\hh_c$
such that $\C{S}_n$ has filtration degree zero,
and any element in $\h\,,\,\h^*$  has filtration degree one, see [BEG].
Then for each $N=1,2,3,\ldots,$ 
the set $Z_N\subset\C$ of all the `$c$' 
for which $P\in \C+[F_N\hh_c,F_N\hh_c]$) is semialgebraic.
Therefore,  either the set  $Z_N$ is
finite or its complement is finite. 
If all sets $Z_N$ are finite then for all `$c$' but countably many, 
$P\notin \C +[\hh_c,\hh_c]$, which contradicts to the fact that 
$P\in \C+[\hh_c,\hh_c]$ for $c\in \C\ssminus T$. Thus, for some 
$N$ the complement to $Z_N$ is finite. Moreover, since 
the set of those `$c$' for which $1\in [F_N\hh_c,F_N\hh_c]$ is finite, 
we see that for all but finitely many $c$, 
we have $P\in \C\oplus [F_N\hh_c,F_N\hh_c]$. 
For such `$c$', there exists a unique complex number $f_P(c)\in \C$ 
(depending on `$c$' as a rational function, of course)
such that $P-f_P(c)\cdot 1_{_{\hh_c}}\in [\hh_c,\hh_c]$. The lemma is proved. 
\end{proof}

\begin{proposition}\label{tracefor}
 Fix $s\in W$. For all $c\in \C$ but countably many, in $\C((z))$ 
we have:
\begin{equation}
\label{trace} 
\Tr(s\cd e^{{z}\cdot \bold h})=\frac{g_{s,c}({z})}{(nc)^{n-1}}\cd \Tr(1), 
\quad\text{where}\quad
g_{s,c}({z})=e^{\frac{(1-nc)(n-1)}{2}z}\cdot
 \frac{\det(1-e^{nc{z}}\cd s)}{\det(1-e^{z}\cd  s)}.
\end{equation}
Moreover, for any $N=1,2,\ldots,$ there exists a finite subset $Y_N\subset \C$ such that 
equation (\ref{trace}) holds modulo ${z}^N$ for all $c\notin Y_N$.  
\end{proposition}

\begin{proof} 
Let $P=s\cdot\bold h^k/k!$. By Lemma \ref{rati}, 
there exists a rational function $f_P$ and 
a finite set $X_P$ such that $\Tr(P)=f_P(c)\cdot\Tr(1)$,
$c\notin X_P$. In particular, we have infinitely many 
$c=r/n$, $r>0$, $(r,n)=1$, which do not belong to $X_P$. 
For such $r$, by Theorem \ref{char}, 
$f_P(c)$ equals the $k$-th term of the series $g_{s,c}({z})/(nc)^{n-1}$
(indeed, it suffices to evaluate the trace in the corresponding 
finite dimensional representation). Since $f_P$ is a rational function, 
this equality holds identically. This implies the Theorem.  
\end{proof}

\begin{remark} Recall the notation
$\fix(s)$ for the dimension of the fixed point set of an element
$s\in S_n$ in
$\h$.
A computation similar to \eqref{t=1} yields $g_{s,c}(0)=(n\cd
c)^{\fix(s)}\,.$
$\lozenge$
\end{remark}

\subsection{Proof of Theorem \ref{classif}}

It is sufficient to prove the result for $c>0$, since 
there is an isomorphism $\hh_c\simeq \hh_{-c}$.

According to \cite[\S3]{BEG}, 
the set of those $c>0$ for which there are 
finite dimensional representations of $\hh_c$  
is of the form ${\sf{S}}+\Bbb Z_+$, where ${\sf{S}}$ is a subset of 
$\lbrace{1/n,...,1-1/n\rbrace}$. Fix $r_0/n\in {\sf{S}}$. 
Then $\hh_c$ has finite dimensional representations
for all $c$ of the form $c=k+(r_0/n)$, $k=0,1,2...$

Take `$k$' to be sufficiently large
and let  $V$ be a finite dimensional representation
of $\hh_{k+(r_0/n)}$. It is clear that $V$ belongs to the category $\oo(\hh_c)$.
Hence,  formula \eqref{char_L}  yields:
\begin{equation*}\label{polynom}
\Tr_{_V}(s\cd  t^{\bold h})=t^{\frac{n-1}{2}}\cd 
\frac{Q_s(t^c)}{\det(1-t\cd s)}\,,\quad\text{where}\quad
Q_s(t)=\sum_{\tau\in\irrep(S_n)}\; a_\tau\cdot\Tr(s,\tau)\cdot 
t^{-{\mathtt{p}}(\tau)}\,.
\end{equation*}  

Further, for each integer $j\geq 0$, we
 consider the finite dimensional representation 
$V_j:=\Bbb S_{c+j-1}\ccirc...\ccirc\Bbb S_c(V)$ of $\hh_{c+j}$. 
Since $c$ is large, for all $i\ge 0,$ the functors $\Bbb S_{c+i}$ are  equivalences
of categories, due to Proposition \ref{sufflarge}.
Hence, the character of $V_j$ is given by the same linear 
combination of characters of standard modules, that is
 given by the formula  
\begin{equation}\label{pole}
\Tr_{_{V_j}}(s\cd  e^{{z} \bold h})=e^{{z}\cdot\frac{n-1}{2}}\cdot
\frac{Q_s(e^{(c+j){z}})}{\det(1-e^{z}\cd  s)}\quad,\quad
\forall s\in S_n\,,\,j=0,1,2,\ldots .
\end{equation}  

On the other hand, by Proposition \ref{tracefor}, for any $N$, 
we have 
\begin{equation}\label{pole2}
\Tr_{_{V_j}}(s\cd  e^{{z}\bold h})=
\dim(V_j)\cd \frac{g_{s,c+j}({z})}{\bigl(n\cdot(c+j)\bigr)^{n-1}}\,
\bigl({\rm mod }\;{z}^N\bigr),
\end{equation} 
for all but finitely many $j$. 
We note that $g_{s,c+j}(0)=
(c+j)^{\fix(s)}$, by  Remark at the end of
\S\ref{trace_sub}.
Observe further that the RHS of formula \eqref{pole}
has a pole at $z=0$ of order $\fix(s)$. Thus, comparing
formulas \eqref{pole} and \eqref{pole2},  we deduce
$$
e^{{z}\cdot\frac{n-1}{2}}\cdot
\frac{Q_s(e^{(c+j){z}})}{\det(1-e^{z}\cd  s)}=
K\cd g_{s,c+j}({z})\,\bigl({\rm mod }\;{z}^N\bigr)\quad
\text{holds for all but finitely many}\enspace j,
$$ 
 where $K$ is a constant
independent of $z\,,\,j,$ and of $s\in S_n$. Now,
since $Q_s$ is a polynomial of fixed degree, and the integer $N$ is
arbitrarily large, 
this implies that $Q_s(t)=K\cd t^{-n(n-1)/2}\cd \det(1-t^n\cdot s)$. 
We conclude that for all $c=r/n$, $r=r_0+kn$, one has 
\begin{equation}\label{char1}
\Tr_{_V}(s\cd  t^\bold h)=K\cd t^{(1-r)(n-1)/2}\cd \frac{\det(1-t^r\cd
s)}{\det(1-t\cd s)}
\quad,\quad
\forall s\in S_n.
\end{equation}

Next, we claim that $(r,n)=1$. To show this, we prove that 
if $(r,n)=d>1$, then the function on the right hand side 
of (\ref{char1}) has poles in $\C^*$. Indeed, let $s\in S_n$ be the
cycle
of order $n$. Then we get 
$$
\Tr_{_V}(s\cd  t^{\bold h})=K\cd t^{(1-r)(n-1)/2}\cd 
\frac{(1-t^{rn})(1-t)}{(1-t^r)(1-t^n)}.
$$
This function has a pole at $t=e^{2\pi i/d}$, which
proves statement (i) of the theorem.

Statement (ii) was proved in the course of the  proof of 
Theorem \ref{reso} in \S3. Namely, it was
shown there that  for all $\tau\ne \triv$ the irreducible representation 
$L(\tau)$ is infinite dimensional.\qed

\subsection{Proof of Theorem \ref{perv}}
For each $i=0,1,\ldots,n-1$ set
$L_i=L(\wedge^i\h)$ and $M_i= M(\wedge^i\h)$.
Let $P_i$ denote the indecomposable
projective cover of $L_i$ in the category $\oo(\hh_c)$, 
which exists by [Gu].
Then the BGG-reciprocity proved in  [Gu], shows that
each $P_i$ has a two-step standard flag with submodule
$M_{i+1}\subset P_i$ and such that $P_i/M_{i+1}\simeq
M_i$ (where $M_n:=0$). Put $P=\oplus_{i=0}^{n-1}\,P_i$.

By general principles, the block of the finite-dimensional represention
in  $\oo(\hh_c)$ is equivalent to the category of finite-dimensional 
representions of the associative algebra
$\A=\Hom_{_{\oo(\hh_c)}}(P,P)$, viewed as a
bigraded algebra $\A=\oplus_{i,j\in\{0,...,n-1\}}\, \A_{ij},$
where $\A_{ij}=\Hom_{_{\oo(\hh_c)}}(P_i,P_j)$.
Furthermore, from the analysis of morphisms in
the category of \linebreak $\HH_q$-modules, that we have carried out earlier
in this section,  one deduces the following

\noindent
{\bf Claim.\;} {\it The algebra $\A$ is isomorphic to
$\C{\mathcal{Q}}/J$,
a
quotient of the {\sl Path algebra} of the following quiver
${\mathcal{Q}}$:
\begin{equation}\label{quiv1}
\xymatrix{
{\underset{0}{\bullet}}\ar@/^/[r]^<>(.5){f_0} &
{\underset{1}{\bullet}}\ar@/^/[r]^<>(.5){f_1}\ar@/^/[l]^<>(.5){g_0} &
{\underset{2}{\bullet}}\ar@/^/[l]^<>(.5){g_1}
\ar@/^/[r]^<>(.5){f_2}&
\;\ldots\;
\ar@/^/[l]^<>(.5){g_2}\ar@/^/[r]^<>(.5){f_{n-3}}
&\ar@/^/[l]^<>(.5){g_{n-3}}
{\underset{n-2}{\bullet}}\ar@/^/[r]^<>(.5){f_{n-2}}&
{\underset{n-1}{\bullet}}\ar@/^/[l]^<>(.5){g_{n-2}},
}
\end{equation}
by the two-sided ideal $J\subset \C{\mathcal{Q}}$ generated
by relations:}
$$
 g_0\ccirc f_0=0\enspace,\enspace
f_{i+1}\ccirc f_{i}=0\enspace,\enspace
g_{i}\ccirc g_{i+1}=0\enspace,\enspace
f_i\ccirc g_i = g_{i+1}\ccirc f_{i+1}\enspace,\enspace
\forall\,i=0,\ldots,n-3\,.
$$

On the other hand, it is well-known that the category
$\perv({\mathbb{P}}^{n-1})$
is also equivalent to the category of finite-dimensional 
representions of the algebra $\C{\mathcal{Q}}/J$.\qed\medskip

This proves Theorem \ref{perv}(i). 
Part (ii) of this Theorem follows from Corollary \ref{length2}, 
and part (iii) from the fact that $M(\tau)$ is simple unless $\tau=
\wedge^i\h$
for some $i$.

\subsection{Proof of Theorems \ref{reso1}\,,\,\ref{BD1}(i), and  Theorem
\ref{char_ehe}(i)-(ii)}

Proposition \ref{sufflarge} and the results of \S2
implies Theorem \ref{reso1}. Indeed, it is easy to see that 
the shift functors for $c=k+\frac{1}{h}$ map standard modules to standard modules. 
Therefore, applying the shift functors $k$ times to the BGG resolution 
at $c=1/h$, we get the BGG resolution for $c=k+\frac{1}{h}$. 
This also implies 
Theorem \ref{BD1}(i), since a representation having the BGG resolution 
is necessarily finite dimensional. Finally, at this point we obtain 
Theorem \ref{char_ehe}(i),(ii), since we have shown that for positive $c=r/n$, 
$(r,n)=1$, the functor $F_c: \oo(\hh_c)\too \oo(\ehe)\,,\,
V \mapsto \e\cdot V$  is an equivalence of categories.

\section{Results for types ${\mathbf{B}=\mathbf{C}}$, and ${\mathbf{D}}$}\label{BDres}
\setcounter{equation}{0}

\subsection{Type ${\mathbf{B}}_n$}

Let $W$ be the Weyl group of type ${\mathbf{B}}_n$, viewed as a subgroup in
$GL(\C^n)$.
 Write $\{e_i\}_{i=1,\ldots,n}$ for the standard basis in $\C^n$.
Then $c=(c_1,c_2)$, where $c_1$
corresponds to roots of the form $\pm e_i\pm e_j$, and $c_2$ corresponds
to roots of the form
to $\pm e_i$.
The Coxeter number of $W$ is $h=2n$.

\begin{theorem} \label{classifi} 
\vi Let $k$ be a nonnegative integer. 
Then for any $c$ such that $c_1(n-1)+c_2=(2k+1)/2$, 
there exists a lowest weight finite dimensional 
representation $N(\triv)$ of $\hh_c$ with lowest 
weight $\triv$ and character given by 
\begin{equation}\label{charfor}
\Tr|_{N(\triv)}(s\cd t^{\bold h})=t^{-kn}\cdot 
 \frac{\det(1-t^{2k+1}\cd s)}{\det(1-t\cd s)}. 
\end{equation}

\vii For all but countably many solutions $c$ of $c_1(n-1)+c_2=(2k+1)/2$,
the $\hh_c$-module $N(\triv)$ is irreducible, i.e., isomorphic to $L(\triv)$.
\end{theorem}

\subsection{Proof of Theorem \ref{classifi}}

Let $u$ be a variable, and
$c_k(u)=(u,\frac{2k+1}{2}-(n-1)u)\in\C[u]^{\oplus 2}$. 
Then $\hh_{c_k(u)}$ is a (flat) algebra over $\C[u]$. 
It follows from Proposition \ref{coxeter}
that there exists a unique 1-dimensional  $\hh_{c_0(u)}$-module
$V_0$ that restricts to the trivial representation of $W$.
We have an isomorphism $V_0=\C[u]$ of $\C[u]$-modules. 

We use the Shift functor
${\Bbb S}_{c,\varepsilon}$, where the character
$\varepsilon: W \to\C^\times$ is given by:
$\varepsilon(s_{e_i-e_j})=1,\varepsilon(s_{e_i})=-1$.
Define an $\hh_{c_k(u)}$-module $V_k:=
{\Bbb S}_{c+(k-1){\mathbf 1}_\varepsilon,\varepsilon}
\ccirc\ldots\ccirc{\Bbb S}_{c,\varepsilon}(V_0)$.

\begin{proposition} \label{coherent} 
$V_k$ is a coherent $\C[u]$-module, whose 
generic fiber is isomorphic to $L(\triv,c)$, and has 
character (\ref{charfor}).
\end{proposition}

\begin{proof}
The fact that $V_k$ is coherent is clear, since 
$\hh_c\e_\varepsilon$ is finite as a right 
$\e_\varepsilon\hh_c\e_\varepsilon$-module.  
We will prove the remaining claims by induction in $k$. 
For $k=0$, the statement is clear, so 
assume it is known for $k=m$ and let us prove it for 
$k=m+1$. 

First of all, $V_{m+1}$ is generically nonzero, 
since $\e V_m$ is generically nonzero (by the character formula for $V_m$). 

On the other hand, let us look at the fiber of $V_{m+1}$ at $u=0$. 
In this case the conclusion of both parts of 
Theorem \ref{classifi} is obvious, since 
$\hh_c=\hh_{0,c_2}=\C S_n\ltimes \hh_{c_2}({\mathbf{A}}_1)^{\otimes n}$. Thus, 
$V_{m+1}|_{u=0}$ is $L(\triv,(0,c_2))$ and its character is given by
formula
(\ref{charfor}).

This implies that $V_{m+1}$ is torsion free at $u=0$.
Indeed, since 
$V_{m+1}$ is generically nonzero, 
its torsion at $u=0$ is a representation of $\hh_{0,c_2}$ 
of smaller dimension than that of $L(\triv,(0,c_2))$, which means 
(by Theorem \ref{classifi} for $c_1=0$) that this torsion must be zero. 

Thus, the fiber of $V_{m+1}$ at generic point $u$ has the character
equal to that of\linebreak
$L(\triv,c_{m+1}(0))$. 
It follows that
 $L(\triv,c_{m+1}(u))$ is contained as a constituent 
in this fiber, and hence has  
the dimension at most that of $L(\triv,c_{m+1}(0))$. 
On the other hand, the module $L(\triv,c_{m+1}(u))$
is the quotient of the standard module $M(\triv,c_{m+1}(u))$
by the radical of  Jantzen-type contravariant form 
in $M(\triv,c_m(u))$. Therefore,
the dimension of $L(\triv,c_{m+1}(u))$ for generic $u$ is at
least the dimension of $L(\triv,c_{m+1}(0))$. 
Hence, these dimensions are equal, and in particular the generic fiber is 
irreducible. 
 
To summarize, we have shown that for generic $u$, 
the fiber at $u$ is an irreducible representation with the character of 
$L(\triv,c_{m+1}(0))$. Thus the induction step has been established,
and the Proposition is proved. 
\end{proof}

Thus, we have shown that $L(\triv,c_k(u))$ is finite dimensional for all 
values of $u$, and generically has the character (\ref{charfor}).
Let $K$ be the kernel of the  contravariant form
in $M(\triv,c_k(u))$,
considered as a $\C[u]$-bilinear form on a free  $\C[u]$-module $M(\triv,c_k(u))$
(thus,  $K$ is itself a free  $\C[u]$-module).
Let $N(\triv,c_k(u))=M(\triv,c_k(u))/K$.
This is a $\C[u]$-free lowest weight $\hh_{c_k(u)}$-module, with character 
(\ref{charfor}). Taking the fibers of this module, we 
get the representations whose existence is claimed 
in Theorem \ref{classifi}(i).

Statement (ii) of the Theorem now  follows immediately from 
the $u=0$ case 
by a deformation argument. Theorem \ref{classifi} is proved.

\begin{remark} In fact, for each $k$ the set of points $(c_1,c_2)$ where 
$N(\triv)\ne L(\triv)$ is finite and consists of rational points. 
The latter follows from a result proved in [DJO], [DO], saying that
the determinant of the contravariant form on weight subspaces of
$M(\triv)$ is a product of linear functions of $c$ with rational coefficients.  
$\lozenge$ 
\end{remark}

\begin{proposition}\label{irre}
If $c=(2k+1)/h$, where $(2k+1,h)=1$, then $N(\triv)$ is irreducible. 
\end{proposition}

\begin{proof}
The character formula \eqref{char_L2} implies that there exist integers
$a_\tau$ such that
$$
\Tr_{_{L(\triv)}}(t^{\bold h})=t^{n/2}\cdot Q(t)/(1-t)^n\quad
\text{where}\quad
Q(t)=\sum_{\tau\in\irrep(W)}\, a_\tau\cdot \dim\tau\cdot t^{-(2k+1)\cdot
{\mathtt{p}}(\tau)/h},
$$
The value of this function at $t=1$, that is 
 the dimension of $L(\triv)$, is given by the expression
 $\frac{(-1)^n}{n!}\cdot (t\frac{d}{dt})^nQ|_{t=1}
=
\sum_{\tau\in\irrep(W)}\,
 a_\tau\cdot (\dim\tau)\cdot [(2k+1)\cdot {\mathtt{p}}(\tau)/h]^n$. 
Since ${\mathtt{p}}(\tau)$ is  an integer, the RHS is clearly
 divisible by $(2k+1)^n$, the dimension of $N(\triv)$. 
Thus, $L(\triv)=N(\triv)$, and we are done. 
\end{proof}

In particular, Theorem \ref{BD1}(ii) is proved. 

\begin{remark}
The same argument as in the proof of 
Proposition \ref{irre} can be used to prove 
the following more general fact. Suppose that $c=(c_1,c_2)$ is a solution of 
the equation $(n-1)\cdot c_1+c_2=(2k+1)/2$, such that $c_1,c_2$ are rational numbers
whose numerators are divisible by $2k+1$ 
(when the numbers are written as irreducible fractions). Then $N(\triv)=L(\triv)$.
$\lozenge$
\end{remark}

\begin{remark} Similar arguments to the ones used in this subsection 
can be employed to study finite dimensional representations 
of $\hh_c({\mathbf{G}}_2)$, using that if $c=(0,c_2)$ then 
$\hh_c({\mathbf{G}}_2)=\Bbb C[\Bbb Z/2\Bbb Z]\ltimes 
H_{c_2}({\mathbf{A}}_2)$. Specifically, 
for any solution $c=(c_1,c_2)$ of the equation $3c_1+3c_2=3k+1$, where $k$ is a 
nonnegative integer,
one can show the existence 
of a finite dimensional lowest weight representation 
$N(\triv)$ with character given by the formula of
Theorem \ref{char}. This is the case
 in particular if $c$ is the constant equal to
$k/2+1/6$. Note that for even $k$ this follows (since $h=6$) from 
Theorem \ref{BD1}(i). Furthermore, in the latter case one has
$N(\triv)=L(\triv)$.   
$\lozenge$
\end{remark}

\subsection{Type ${\mathbf{D}}_n$} 

Now let $W$ be of type ${\mathbf{D}}_n$. 
In this case the Coxeter number is $h=2(n-1)$. 

\begin{theorem}\label{classifi1} 
\vi If $c=(2k+1)/h$ where $k$ is a nonnegative integer, then there exists a
finite dimensional lowest weight module 
$N(\triv)$ over $\hh_c$ with character 
 given by (\ref{charfor}).

\vii If $(2k+1,h)=1$ then $N(\triv)$ is irreducible (i.e. is isomorphic to $L(\triv)$).  
\end{theorem}

\begin{proof}
Recall that $\hh_{c,0}({\mathbf{B}}_n)=\C[\Bbb Z/2]\ltimes \hh_c({\mathbf{D}}_n)$. 
By Theorem \ref{classifi}, for $c=r/h$, 
there is a 
lowest weight representation $N(\triv)$ of $\hh_{c,0}({\mathbf{B}}_n)$,
which has character given by (\ref{charfor}). 
Restricting this representation to $\hh_c({\mathbf{D}}_n)$, we obtain a required 
representation of $\hh_c({\mathbf{D}}_n)$.

The proof of its irreducibility for $(2k+1,n)=1$ is 
analogous to the proof of Proposition \ref{irre}. The theorem is proved.
\end{proof} 

In particular, Theorem \ref{BD1}(iii) is proved. 

\subsection{An example}

An example of the algebra $\hh_c({\mathbf{D}}_4)$ 
shows, as we will now see, that in the simplest
non-relatively prime case $c=1/2$ the module $N(\triv)$, constructed
in the Theorem above, is {\it not} irreducible, and actually contains a submodule
whose lowest weight is the reflection representation tensored with a character. 

In more detail, we start with the algebra $\hh_c({\mathbf{B}}_4)$ with
parameter $c=(c_1,c_2)$ such that $6c_1+2c_2=3.$
Then by Theorem \ref{classifi},
$N(\triv)$ is an 81-dimensional lowest weight representation
of  $\hh_c({\mathbf{B}}_4)$,
with character $(t^{-1}+1+t)^4.$
Generically, this representation is irreducible.
Moreover, for  $c_1=0$ hence generically,
we have an $\hh_c$-module map $\psi: M(\h)\to M(\triv)$,
and $N(\triv)= M(\triv)/{\mathtt{Image}}(\psi)$,
where $
{\mathtt{Image}}(\psi)$ is generated by a copy of $\h$ sitting in degree 3 in $M(\triv)=
\Bbb C[\h].$
There are two copies of $\h$ in $
M(\triv)=\C[\h]$, obtained by taking the partial derivatives of
the two invariant polynomials in degree $4.$
There is a copy of $\h\subset \C[\h]=M(\triv)$ consisting of singular vectors, i.e. those
annihilated by the elements of $\h\subset \hh_c$ acting in
$\C[\h]$ via the  Dunkl operators. This copy 
is  spanned by partial derivatives of some polynomial $Q$ of degree 4,
invariant under $W= W({\mathbf{B_4}})$. Such a polynomial has the form
$a\cdot(\sum_{i=1}^3\, x_i^2)^2-\sum_{i=1}^3\,
x_i^4.$ An explicit calculation with Dunkl operators shows that
$a=c_1$. So, generically $L(\triv)$ is the local ring
of the isolated critical point of $Q$ at the origin.

However, sometimes  the origin is not an isolated critical point, in which case 
the local ring is not finite dimensional, and 
hence cannot be isomorphic to $L(\triv)$.
This happens, for instance, if $c_1=1/p\,,\, p =1,..,4;$ in particular,
 for $c_1=1/2$
(in which case $c_2=0$, so we are essentially in the ${\mathbf{D}}_4$ case).
In this case, we have additional polynomials in degree $3$,
killed by Dunkl operators. 
An explicit calculation shows
that this additional space of singular vectors in $\C[\h]=M(\triv)$
is isomorphic, as a representation of $W$, to $\h\otimes \varepsilon$
where $\varepsilon$ is the character
such that $\varepsilon(s_{e_i-e_j})=1$
and $\varepsilon(s_{e_i})=-1$.  
More specifically, a basis of this space is formed by the polynomials
$f_i=\prod_{j\ne i}\,x_j$. 
These polynomials are definitely nonzero in $N(\triv)$, which shows that 
$N(\triv)\ne L(\triv)$. In fact, it is easy to see 
that we have an exact sequence
$$
0\to L(\h\otimes \varepsilon)\to N(\triv)\to L(\triv)\to 0.
$$

\section{Connections to other topics}
\setcounter{equation}{0}
\subsection{Relation to the affine flag variety.}\label{Chere}

Cherednik introduced  an associative
 $\C$-algebra depending on two
complex parameters $q,t \in \C^*$, called the
{\it double affine Hecke algebra}, see e.g.   [Ch1]
The double affine Hecke algebra
${\mathbf{H}}_{q,t}$ has generators: $T_i\,,\, X_i^{\pm 1}\,,\, Y_i^{\pm
1}\,,\,
i=1,\ldots,\ell,$ with certain defining relations, 
analogous to those in the affine Hecke algebra
(the latter is a subalgebra in ${\mathbf{H}}_{q,t}$
generated by the $T_i$'s and $X_i^{\pm 1}$'s).
Cherednik showed 
 that the elements
$\{T_w\cdot X_i^a\cdot Y_j^b\;|\;w\in W\,;\,
i,j=1,\ldots,\ell\,;\,a,b\in \Z\},$ form
a Poincar\'e-Birkhoff-Witt type $\C$-basis in ${\mathbf{H}}_{q,t}$.

The rational Cherednik algebra $\hh_c$ may be thought of,
see [EG], as the following
two-step degeneration of the double affine Hecke algebra:
$${\mathbf{H}}_{q,t}\;=\;{\mathbf{H}}^{^{\,_{\text{elliptic}}}}_{q,t}
\enspace\rightsquigarrow \enspace
{\mathbf{H}}^{^{\,_{\text{trigon}}}}_c\enspace\rightsquigarrow \enspace
{\mathbf{H}}^{^{\,_{\text{rational}}}}_c\;=\;\hh_c\,.
$$
Otherwise put, the algebra  ${\mathbf{H}}_{q,t}$ 
can be regarded as a formal $\C[[\epsilon]]$-deformation 
of $\hh_c$ upon the substitution:
 $t=e^{{\epsilon}^2c}$, $X_i=e^{{\epsilon}\cdot x_i}$, 
$Y_i=e^{{\epsilon}\cdot y_i}$.
 In \cite[p.65]{Ch2}, it is observed that 
this deformation is in effect trivial (as had been conjectured by 
the second author).

 Before our paper 
was written, Cherednik announced \cite[p.64]{Ch2}
a classification of finite dimensional 
representations of the double affine Hecke algebra 
${\mathbf{H}}_{q,t}$ of type ${\mathbf{A}}$.
 As was pointed out in [Ch2], this allows one 
to classify finite dimensional representations
of $\hh_c(S_n)$ (i.e. our Theorem 1.2), 
and obtain the dimension formula for them
using the corresponding results for  the double affine Hecke algebra.

Let us explain, for example, why the trigonometric Cherednik
algebra ${\mathbf{H}}^{^{\,_{\text{trigon}}}}$ 
is trivial as a formal deformation of the rational
Cherednik algebra $\hh_c $
(the argument is due to Cherednik). Recall that the
(formal) trigonometric Cherednik algebra
${\mathbf{H}}^{^{\,_{\text{trigon}}}}(\epsilon)$, 
can be realized,
via the Dunkl representation, as
a complete topological $\C[[\epsilon]]$-algebra (topologically)
generated by
$W$, $\C[\h]$, and trigonometric Dunkl operators
$$
T^{^{\,_{\text{trigon}}}}_y(\epsilon)=\partial_y-
\sum\nolimits_{\alpha\in R_+}c_\alpha\cd 
\frac{\epsilon\cd\langle\alpha,y\rangle}{1-e^{-\epsilon\cdot\alpha}}\cd(1-s_\alpha)
+\frac{\epsilon}{2}\sum\nolimits_{\alpha\in R_+}\,
c_\alpha\cdot\langle\alpha,y\rangle\;,\quad
y\in \h$$
($\epsilon$ is a deformation parameter).
Since the function ${1\over x}-
{\epsilon\over (1-e^{-\epsilon\cdot x})}$ is regular at
$\epsilon=0$,
we find
$$T^{^{\,_{\text{trigon}}}}_y(\epsilon)=\partial_y-\sum\nolimits_{\alpha\in R_+}c_\alpha\cd 
\frac{\langle\alpha,y\rangle}{\alpha}\cd(1-s_\alpha)
+ \epsilon\cdot F=
T_y^{^{\,_{\text{rational}}}}+ \epsilon\cdot F\,,
$$
where
 $T^{^{\,_{\text{rational}}}}_y$ is the {\it rational} Dunkl operator
and $F\in 
\C{W}\ltimes \C[\h][[\epsilon]]$.
It follows readily that there is
a topological  $\C[[\epsilon]]$-algebra isomorphism:
${\mathbf{H}}^{^{\,_{\text{trigon}}}}(\epsilon) \simeq
\C[[\epsilon]]{\hat\otimes}
\hh_c.$
We see that, if for some $c,$ the algebra  $\hh_c,$
 has representation in a $d$-dimensional
$\C$-vector space ($d<\infty$), then  the $\C[[\epsilon]]$-algebra
${\mathbf{H}}^{^{\,_{\text{trigon}}}}(\epsilon) \simeq
\C[[\epsilon]]{\hat\otimes}
\hh_c$ has a $\C[[\epsilon]]$-linear representation in 
 $\C[[\epsilon]]^d$. Furthermore, we observe that  all the operators in
$\C[[\epsilon]]^d$
corresponding to
the generators $x_i\,,\,i=1,\ldots,\ell,$ of $\hh_c$ are necessarily
nilpotent, i.e. the infinite series involved
are actually polynomials in $\epsilon$.
 Hence, it makes sense to specialize the formal variable
$\epsilon$ to the complex number: $\epsilon=1,$ to
obtain
a well-defined action in $\C^d$ of the original $\C$-algebra
${\mathbf{H}}^{^{\,_{\text{trigon}}}}_c$. This way, one proves
\begin{proposition}[Cherednik]\label{formal}
There is a natural  equivalence of the categories
of finite dimensional $\hh_c$- and
${\mathbf{H}}^{^{\,_{\text{trigon}}}}_c$-modules
that preserves the $\C$-dimension of the modules.\qed
\end{proposition}

Motivated by  the well-known geometric interpretation
of the affine Hecke algebra (see [CG, ch. 7] for a review),
E. Vasserot [Va] has produced a similar construction
for the algebra ${\mathbf{H}}_{q,t}$. In more detail,
let $R$ be the root system of a  complex  semisimple 
 Lie
algebra  $\bar\g$,
let $\bg={\bar\g}((z))$ be the corresponding loop Lie
algebra, and $\GG$ the corresponding Kac-Moody group.
One has the Springer resolution $\tg \to \bg$ and
the Steinberg variety $\zz=\tg\times_\bg\tg$. All the
(infinite dimensional) varieties $\bg, \tg,$ and $\zz$
acquire a natural $\GG$-action by conjugation,
and a natural $\C^*$-action by dilations,
thus become $\GG\times\C^*$-varieties.
 Vasserot [Va] succeeded in
defining an equivariant $K$-group $K^{\GG\times\C^*}(\zz)$ with
a convolution-type (non-commutative) algebra structure,
and proved (roughly speaking) an algebra
isomorphism ${\mathbf{H}} \simeq K^{\GG\times\C^*}(\zz)$,
where ${\mathbf{H}}$ is the double affine Hecke algebra
associated to the dual root system $R^\vee$.

It is very likely that the argument in [Va] can be used
to prove a similar  algebra
isomorphism ${\mathbf{H}}^{^{\,_{\text{trigon}}}}
 \simeq H^{\GG\times\C^*}(\zz)$, where $H^{\GG\times\C^*}(-)$ stands for 
$\GG\times\C^*$-equivariant Borel-Moore homology. We note at this point
that {\it no} analogous interpretation of the {\it rational}
Cherednik
algebra ${\mathbf{H}}^{^{\,_{\text{rational}}}}_c=\hh_c$
is available at present.

 Using the isomorphism
${\mathbf{H}} \simeq K^{\GG\times\C^*}(\zz)$, Vasserot obtained
 a complete classification of
simple `bounded below' ${\mathbf{H}}$-modules in geometric
terms. A similar classification of  
simple `bounded below'
${\mathbf{H}}^{^{\,_{\text{trigon}}}}$-modules,
based on the isomorphism ${\mathbf{H}}^{^{\,_{\text{trigon}}}}
 \simeq H^{\GG\times\C^*}(\zz)$ should be possible.
Although  classifications of that type do {\it not} 
allow, in general, to distinguish finite dimensional
representations among all  `bounded below' representations,
they do allow to construct some of them. In this manner,
some simple  finite dimensional ${\mathbf{H}}$-modules 
have been constructed in \cite[\S9.3]{Va}. It seems certain that
their trigonometric analogues exist as well. 
By Proposition \ref{formal},
 finite dimensional
${\mathbf{H}}^{^{\,_{\text{trigon}}}}_c$-modules 
arising in that way
may also be seen as $\hh_c$-modules. A comparison
with \cite[\S9.3]{Va} shows that, in type $\mathbf{A}$,
the geometric construction of {\it loc. cit.} produces in effect
 all simple  finite dimensional $\hh_c$-modules.
In other types, this is definitely not true. The complexity
of the situation for types $\mathbf{B}, \mathbf{D},$
observed both in the present paper and in [Va], 
suggests that some kind of `$L$-packet phenomenon'
resulting from the failure   of $\GG$-orbits to be
simply-connected is involved.

Let $\bb$ denote the flag variety for $\GG$, the so-called
affine flag manifold.
Recall that simple ${\mathbf{H}}$-modules have been realized
in [Va] as subquotients of $H^*(\bb_x^s)$, the cohomology
of certain fixed point sets $\bb_x^s$ in $\bb$. In the case
of a {\it nil-elliptic} element $x\in\bg$ the corresponding
set $\bb_x^s$ turns out to be a (finite dimensional) smooth projective
variety, see [KL]; moreover, the  total cohomology $H^*(\bb_x^s)$
 turns out to be a simple  ${\mathbf{H}}$-module.
A trigonometric analogue of this 
should provide a realization of finite dimentional simple
${\mathbf{H}}^{^{\,_{\text{trigon}}}}_c$-modules,
hence of $\hh_c$-modules,  as the total cohomology groups
of appropriate smooth projective subvarieties in $\bb$.
The  Poincar\'e duality for such varieties
would
explain the
Gorenstein property in Proposition \ref{gor}.

Let  $x\in \bg$ be  a nil-elliptic element.
We note that the computation of the Poincar\'e polynomials
for varieties $\bb_x^s$ 
(carried out in [LS] in type $\mathbf{A}$ and in [So] in general)
produces exactly the same
answer as our formula in Theorem \ref{char}.
In particular, for the Euler characteristic,
the computation   in [LS],
[So] gives: $\chi(\bb_x^s)= r^\ell$,
where the integer $r$ is  related to $x$
as explained e.g. in \cite[\S8]{Va}.
This agrees with our Proposition \ref{traceCorollary}.

Replacing the affine flag manifold
by a loop Grassmannian $\Gr$, one constructs similarly
the `spherical' Steinberg variety $\zz\spher$.
One expects to have algebra isomorphisms
$\e{\mathbf{H}}\e \simeq K^{\GG\times\C^*}(\zz\spher),$
and $\e{\mathbf{H}}^{^{\,_{\text{trigon}}}}\e \simeq 
H^{\GG\times\C^*}(\zz\spher)$. This should give rise to a geometric
construction of simple `bounded below'
$\e{\mathbf{H}}\e$-modules and 
$\e{\mathbf{H}}^{^{\,_{\text{trigon}}}}\e$-modules,
respectively. Our expectations are supported by
the computations of
 the Euler characteristic of certain fixed point varieties
$\Gr_x^s$ carried out in [So]. The result of Sommers reads:
$\chi(\Gr_x^s)= \prod_{i=1}^\ell\,
\frac{d_i+r-1}{d_i},$
which is according to our Theorem \ref{char_ehe} exactly 
the dimension of a simple finite dimensional $\ehe$-module.

\subsection{Relation to the geometry of Hilbert scheme}

  Let ${\sf{Hilb}}^n(\C^2)$ be the {\it Hilbert
scheme} of $n$-points in $\C^2$, that is, a scheme
parametrizing all codimension $n$ ideals in the polynomial ring
$\C[x,y]$
in two variables. It is known that  ${\sf{Hilb}}^n(\C^2)$
is a smooth  connected algebraic variety of dimension $2n$.
There is a natural (ample) determinant line bundle
${\mathcal{L}}$ on ${\sf{Hilb}}^n(\C^2)$, 
the $n$-th wedge power of the tautological bundle of rank $n$.
 The natural action on the plane $\C^2$ of the torus
$\TT=\C^*\times\C^*$, by diagonal matrices,
lifts to a $\TT$-action on ${\sf{Hilb}}^n(\C^2)$,
making ${\mathcal{L}}$  a $\TT$-equivariant line bundle.

There is a natural Hilbert-Chow
morphism $\pi: {\sf{Hilb}}^n(\C^2)\too \sym^n(\C^2)=\C^{2n}/S_n.$
This morphism is proper, and one puts
${\sf{Hilb}}^n_o(\C^2):=\pi^{-1}(0),$
the zero-fiber of $\pi$.
The scheme ${\sf{Hilb}}^n_o(\C^2)$ parametrizes all codimension $n$
ideals in  $\C[x,y]$ set-theoretically concentrated at the origin $0\in
\C^2$;
alternatively,  ${\sf{Hilb}}^n_o(\C^2)$ parametrizes  codimension $n$
ideals in the 
formal power series ring $\C[[x,y]]$. It is known that  ${\sf{Hilb}}^n_o(\C^2)$
is a reduced and irreducible  $\TT$-stable subscheme in ${\sf{Hilb}}^n(\C^2)$
of dimension $(n-1)$.

Giving an algebraic $\TT$-action on a vector space $V$ is equivalent
to giving a bigrading $V= \oplus_{i,j\in\Z}\, V_{ij}$. If all the
bigraded components are finite-dimensional,
we may introduce
the formal character
 $\cht(V):=\sum_{i,j\in\Z}\, q^i\cdot t^j\cdot \dim V_{ij}
\,\in\C[[q,q^{-1},t,t^{-1}]]$.

The results of M. Haiman 
[Ha2],[Ha3] imply the following 2-variable character formula
\begin{equation}\label{Ha1}
\cht\bigl( H^0({\sf{Hilb}}^n_o(\C^2)\,,\,{\mathcal{L}}^{\otimes k})\bigr)
= C^{(k)}_n(q,t)\quad,\quad\forall k=0,1,2,\ldots\;,
\end{equation} 
where the $C^{(k)}_n(q,t)$ are $(q,t)$-analogues of Catalan numbers
introduced in \cite[(1.10)]{Ha2}. In particular, the dimension of the
cohomology space  is given by the usual Catalan number:
\begin{equation}\label{Hadim}
\dim
H^0({\sf{Hilb}}^n_o(\C^2)\,,\,{\mathcal{L}}^{\otimes k})
= \frac{1}{n\cdot k+1}\cd \left({n(k+1)\atop n}\right)
\quad,\quad\forall k=0,1,2,\ldots\;.
\end{equation} 

Now, given $c\in \C$, write $L_c(\triv)=L(\triv)$ (to
keep track of the parameter `$c$') for the simple $\hh_c(S_n)$-module
corresponding to the trivial $S_n$-representation.
We see that the  RHS of formula \eqref{Hadim} equals, by Theorem
\ref{traceCorollary}, the dimension
of the simple  finite dimensional  $\ehe$-module:
\begin{equation}\label{Hadim2}
\dim
H^0({\sf{Hilb}}^n_o(\C^2)\,,\,{\mathcal{L}}^{\otimes k})
= \dim\bigl(\e\cd L_c(\triv)\bigr)\quad\text{for}\quad
c=\frac{1}{n}+k\;,\;k=0,1,2,\ldots\;.
\end{equation} 
We will now provide a heuristic `explanation' of 
 equality \eqref{Hadim2}.

Let $\h=\C^n$. Write $\C[\h\oplus\h]^W$
for the space  of 
diagonal $W$-invariants.
Let $\m=\C[\h\oplus\h]^W_+$ denote the augmentation ideal in $\C[\h\oplus\h]^W$
formed by all
diagonal $W$-invariants without constant term,
and $\C_o :=\C[\h\oplus\h]^W/\m$, the
corresponding
1-dimensional module. 

The Hilbert-Chow morphism gives a canonical isomorphism
$H^0({\sf{Hilb}}^n(\C^2)\,,\,
\OH)=\C[\h\oplus\h]^W$. Thus,
the space $H^0({\sf{Hilb}}^n(\C^2)\,,\,{\mathcal{L}}^{\otimes k})$
acquires a natural $\C[\h\oplus\h]^W$-module structure.
 Furthermore, thanks to \cite[(25)]{Ha3}, we have:
$$H^0({\sf{Hilb}}^n_o(\C^2)\,,\,{\mathcal{L}}^{\otimes k})
= H^0({\sf{Hilb}}^n(\C^2)\,,\,{\mathcal{L}}^{\otimes k})/\m\cdot
H^0({\sf{Hilb}}^n(\C^2)\,,\,{\mathcal{L}}^{\otimes k}).$$
Therefore, the space on the LHS of \eqref{Hadim2} may be
rewritten as follows
\begin{equation}\label{surj}
H^0({\sf{Hilb}}^n_o(\C^2)\,,\,{\mathcal{L}}^{\otimes k})\,\simeq\,
H^0\Bigl({\sf{Hilb}}^n(\C^2)\,,\,\underbrace{{\mathcal{L}}
\otimes_{_{\OH}}
\ldots\otimes_{_{\OH}}{\mathcal{L}}}_{k\;{\footnotesize{{times}}}}\Bigr)
\bigotimes\nolimits_{\C[\h\oplus\h]^W}\,\C_o\,.
\end{equation}

%%%%%%%%
We now turn to the RHS of  \eqref{Hadim2}. 
Given $c=\frac{1}{h}+k,$ introduce the shorthand
 notation $\hh^{(k)} = \hh_{h^{-1}+k}.$
By the definition of {\it Shift functor},
see Lemma \ref{morita} and
Proposition \ref{sufflarge}(ii),
 for $c=\frac{1}{h}+k$ where $k=0, 1,\ldots\,,$  we have
\begin{align}\label{ha4}
L_c(\triv) & = {\mathbb{S}}_{k-1+h^{-1}} \ccirc\ldots 
\ccirc{\mathbb{S}}_{1+h^{-1}}\ccirc{\mathbb{S}}_{h^{-1}}(\C)\\
 & = \hh^{(k)}{\emi}
\bigotimes\nolimits_{\e\hh^{(k-1)}\e}\e\hh^{(k-1)}{\emi} 
\bigotimes\nolimits_{\e\hh^{(k-2)}\e}\ldots\bigotimes\nolimits_{\e\hh^{(1)}\e}
\e\hh^{(1)}{\emi} \bigotimes\nolimits_{\e\hh^{(0)}\e}\C.\nonumber
\end{align}
where we have repeatedly exploited the 
algebra isomorphism: $\emi\hh^{(k+1)}\emi \simeq \e\hh^{(k)}\e$,
see Proposition \ref{isom}. Similarly,
one has:
\begin{equation}\label{ha4e}
\e\cd L_c(\triv)
= \left(\e\hh^{(k)}{\emi}
\bigotimes\nolimits_{\e\hh^{(k-1)}\e}\e\hh^{(k-1)}{\emi} 
\bigotimes\nolimits_{\e\hh^{(k-2)}\e}\ldots\bigotimes\nolimits_{\e\hh^{(1)}\e}
\e\hh^{(1)}{\emi}\right) \bigotimes\nolimits_{\e\hh^{(0)}\e}\C.
\end{equation}

Next,
recall  that the algebra $\hh_c$ has a canonical increasing filtration
such that
the generators $x\in\h^*\,,\,y\in\h$ are assigned filtration degree 1,
and elements of the Weyl group $W$  are assigned filtration degree 0.
The filtration on
$\hh_c$ induces  natural increasing
filtrations on other objects. 
For the corresponding associated graded objects one  has:
\begin{align}\label{grd}
&\grd\hh_c\simeq
\C[\h\oplus\h^*]\#W\quad,\quad
\grd(\ehe)\simeq \C[\h\oplus\h]^W \simeq
H^0\bigl({\sf{Hilb}}^n(\C^2)\,,\,\OH\bigr) \quad\text{and}\nonumber\\
&\hphantom{x}\hspace{120pt}\grd(\e\hh_c\emi) \simeq
\C[\h\oplus\h]^\sign\simeq H^0\bigl({\sf{Hilb}}^n(\C^2)\,,\,{\mathcal{L}}\bigr)\,.
\end{align}
The last two isomorphisms in  \eqref{grd}
suggest to think of \eqref{ha4e} as a `quantum' analogue 
of \eqref{surj}.

To make this analogy more precise, observe first
that the
 canonical filtration on  $\hh_c$ is stable under the adjoint action
$\ad\bh: u\mapsto [\bh, u]=\bh\cdot u - u\cdot\bh,$
hence induces an $\ad\bh$-action on the associated graded algebra.
It is immediate to check that transporting this action under the first
isomorphism in
\eqref{grd} we get: $\ad\bh(x)=x\,,\,\forall
x\in\h^*,$
and $\ad\bh(y)=-y\,,\,\forall y\in\h$.
Thus,   $\grd\hh_c$ may be viewed as a  {\it
bi-graded} algebra $\grd\hh_c=\oplus_{p,q\in\Z}\,\hh^{p,q}$,
where the space  $\hh^{p,q}$ is defined to have total degree
$p+q$ and such that $\ad\bh|_{_{\hh^{p,q}}}=
(p-q)\cdot\id_{_{\hh^{p,q}}}.$
Then, all the isomorphisms in \eqref{grd} become bigraded isomorphisms,
where the bidegrees on $\C[\h\oplus\h^*]\#W$ are given by:
$\deg w=(0,0)\,,\,\deg x=(1,0),$ and $\deg y= (0,1)\,,\,
\forall w\in W\,,\,x\in\h^*\,,\,y\in\h.$

The canonical filtration on
$\hh_c$ induces a natural increasing
filtration  $F_\bullet$ on each term  of the tensor product in the RHS of 
\eqref{ha4}. We introduce an  increasing
filtration on $L_c(\triv)$ to be the induced filtration
\begin{align*} 
&F_\bullet(L_c(\triv)) :=\hspace{100pt} \hphantom{x}\\
&\hphantom{x}\hspace{30pt} F_\bullet\left(\hh^{(k)}{\emi}
\bigotimes\nolimits_{\e\hh^{(k-1)}\e}\e\hh^{(k-1)}{\emi} 
\bigotimes\nolimits_{\e\hh^{(k-2)}\e}\ldots\bigotimes\nolimits_{\e\hh^{(1)}\e}
\e\hh^{(1)}{\emi} \bigotimes\nolimits_{\e\hh^{(0)}\e}\C\right),
\end{align*}
where, for any $a\in\Z$, the corresponding $a$-th term of the
filtration on the RHS above is defined as the image of the $\C$-vector
space
\begin{align*}
\sum_{\{a_1,\ldots,a_k\geq 0\;|\; a_1+\ldots+a_k\leq a\}}\,
F_{a_k}\bigl(\hh^{(k)}{\emi}\bigr)\otimes_{_\C}
F_{a_{k-1}}\bigl(\e\hh^{(k-1)}{\emi}\bigr)\otimes_{_\C}
 \ldots\otimes_{_\C}
F_{a_1}\bigl(\e\hh^{(1)}{\emi}\bigr)\otimes_{_\C}\C\\
\subset \;
\hh^{(k)}{\emi}
\otimes_{_\C}\e\hh^{(k-1)}{\emi} 
\otimes_{_\C}\ldots\otimes_{_\C}
\e\hh^{(1)}{\emi} \otimes_{_\C}\C
%\nonumber\\ 
%&:= \sum_{\{a_1,\ldots,a_k\geq 0\;|\; a_1+\ldots+a_k\leq a\}}\,
%F_{a_k}\bigl(\hh^{(k)}{\emi}\bigr)\otimes
%F_{a_{k-1}}\bigl(\e\hh^{(k-1)}{\emi}\bigr)\otimes
% \ldots\otimes
%F_{a_1}\bigl(\e\hh^{(1)}{\emi}\bigr)\otimes\C.\nonumber
\end{align*}
under the canonical
projection from a tensor product over $\C$ to the corresponding
tensor product over the algebras $\e\hh^{(i)}\e$.

The filtration on $L_c(\triv)$ thus defined
is clearly $\bh$-stable.\footnote{We remark that the lowest weight line
$\triv\subset L_c(\triv)$ does {\it not} usually belong to the
first (the smallest) term of the filtration on  $L_c(\triv)$.}  Therefore, we obtain 
a well-defined $\bh$-action on the associated graded space.
As we have done above in the case of the algebra $\hh_c$ itself, 
this allows  to view $\grd L_c(\triv)$
 as a
bi-graded space
\begin{equation}\label{grd2}
\grd L_c(\triv)=\oplus_{p,q\in\Z}\,L_{p,q}\quad\text{where}\quad
\bh|_{_{L_{p,q}}}
=(p-q)\cdot\id_{_{L_{p,q}}}\,,
\end{equation}
and such that the original grading corresponds to the {\it total}
grading, that is, we have:
 $\grd^aL(\triv)=\oplus_{\{(p,q)\;|\;p+q=a\}}\,L_{p,q}.$

Similar considerations and definitions apply to the $\ehe$-module $\e\cdot
L(\triv)$. In particular, we define an increasing filtration
on $\e\cd L(\triv)$ using the isomorphism of \eqref{ha4e}, and view
$\grd\bigl(\e\cdot L(\triv)\bigr)$  as a
bi-graded space, that is as a $\TT$-module.
Thus, we may define the {\it bigraded character} of $\grd\bigl(\e\cdot  L(\triv)\bigr)$,
a two-variable formal series 
$\chi_{_{\e\cdot L(\triv)}}(q,t)\in \C[q,q^{-1},t,t^{-1}]$.

The following is a considerable refinement of  equation
\eqref{Hadim2}.
\begin{conjecture}\label{quant}
{\it Let $c=\frac{1}{n}+k$. Then there is a canonical $\TT$-module isomorphism
$$
\grd\bigl(\e\cd L_c(\triv)\bigr)\; \simeq\;
H^0({\sf{Hilb}}^n_o(\C^2)\,,\,{\mathcal{L}}^{\otimes k})
\,,\quad\text{for any}\quad k=0, 1,\ldots\,.
$$
In particular, $\chi_{_{\e\cdot L(\triv)}}(q,t)=C^{(k)}_n(q,t)$,
is the $(q,t)$ Catalan number.}
\end{conjecture}

To formulate an analogue of Conjecture \ref{quant} for the $\hh_c$-module
$L_c(\triv)$, we need to recall
(see [Ha3]) that there is an `unusual' tautological
rank $n!$ vector bundle $\rr$ on ${\sf{Hilb}}^n(\C^2)$ whose
fibers afford the regular representation of the group $W=S_n$.
\footnote{Haiman denotes this vector bundle by `$P$', in honor of Procesi.}

\begin{conjecture}\label{quant2}
{\it Let $c=\frac{1}{n}+k+1.$
Then there is a canonical $W$-equivariant $\TT$-module isomorphism}
$$
\grd\bigl(L_c(\triv)\bigr)\; \simeq\;\sign\otimes
H^0({\sf{Hilb}}^n_o(\C^2)\,,\,\rr \otimes {\mathcal{L}}^{\otimes k})
\,,\quad\text{for any}\quad k=0,1, \ldots\,.
$$
\end{conjecture} 
Case $k=0$ of  Conjecture \ref{quant2} is true; it  amounts (modulo [Ha3])
to the main result of Gordon [Go] on diagonal harmonics.
Conjecture \ref{quant} is easy for $k=0$;
Case $k=1$ of Conjecture \ref{quant}
follows from  Conjecture \ref{quant2} for  $k=0$.

\begin{remark}
If true, Conjecture \ref{quant2} would provide a formula, see \cite[(106)]{Ha3}, for the
bigraded multiplicities 
$\,\dis [\grd\bigl(L_c(\triv)\bigr){}:{}\tau]\,,\, \tau \in
\irrep(S_n),$
in terms
involving among other things
the $(q,t)$-Kostka polynomials, whose positivity has been
conjectured by Macdonald and proved by Haiman. $\lozenge$
\end{remark}
\medskip

\begin{remark}
The geometric structures considered in \S6.1 and  \S6.2 should be
related  to each other by a kind of `Langlands duality' for
 rational Cherednik
algebras, similar somewhat to the existing Langlands duality
for affine Hecke algebras, cf. e.g. \cite[Introduction]{CG}. These matters are 
not understood at the moment.
\end{remark}

\footnotesize{
{\bf Y.B.}: Department of Mathematics, Cornell University,
Ithaca, NY 14853-4201, USA;\\
\hphantom{x}\quad\, {\tt berest@math.cornell.edu}

{\bf P.E.}: Department of Mathematics, Rm 2-165, MIT,
77 Mass. Ave, Cambridge, MA 02139;\\
\hphantom{x}\quad\, {\tt etingof@math.mit.edu}

{\bf V.G.}: Department of Mathematics, University of Chicago,
Chicago, IL
60637, USA;\\
\hphantom{x}\quad\, {\tt ginzburg@math.uchicago.edu}}

\end{document}